\documentclass[12pt]{amsart}
\usepackage{amsmath, amsthm,txfonts}
\usepackage{tikz}

\newtheorem{theorem}{Theorem}[section]
\newtheorem{lemma}[theorem]{Lemma}
\newtheorem{proposition}[theorem]{Proposition}
\newtheorem{corollary}[theorem]{Corollary}

\theoremstyle{definition}

\newtheorem{example}[theorem]{Example}
\theoremstyle{remark}
\newtheorem{remark}[theorem]{Remark}
\numberwithin{equation}{section}

\title{Spectra of three-peg Hanoi towers graphs}

\author[Hungar]{Brett Hungar}
\address{Department of Mathematics, Ohio State University, 231 W 18th Ave, Columbus, OH 43210}
\email{bretthungar@gmail.com}

\author[Mograby]{Gamal Mograby}
\address{Department of Mathematics, University of Maryland, College Park, MD 20742-4015}
\email{gmograby@umd.edu} 

\author[Phelps]{Madison Phelps}
\address{Department of Mathematics, Kidder Hall 368, Oregon State University, Corvallis, OR 97331-4605}
\email{phelpmad@oregonstate.edu}

\author[Rogers]{Luke~G. Rogers}
\address{Department of Mathematics, University of Connecticut, Storrs, CT, 065269-1009}
\email{luke.rogers@uconn.edu}

\author[Wheeler]{Jonathan Wheeler}
\address{Department of Mathematics, Northern Arizona University, 801 S. Osborne Dr. PO Box: 5717, Flagstaff, AZ 86011}
\email{jlw697@nau.edu}


\begin{document}

\begin{abstract}
We consider the relationship between the Laplacians on two sequences of planar graphs, one from the theory of self-similar groups and one from analysis on fractals. By establishing a spectral decimation map between these sequences we give an elementary calculation of the spectrum of the former, which was first computed by Grigorchuk and \v{S}uni\'{c}~\cite{GrigorchukSunic}. Our method also gives a full description of the eigenfunctions.
\end{abstract}

\maketitle

\section{Introduction}

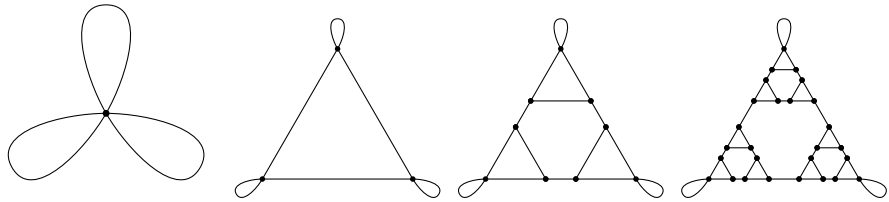
\begin{figure}
\begin{tikzpicture}[scale=0.5]
\fill  (2.5, 1.44338) circle[radius=2.6pt] ;
\draw (2.5,1.44338) to  [out=180,in=120] (0,0) node{} to [out=300, in=240] (2.5,1.44338);
\draw (2.5, 1.44338) to [out=0,in=60]  (5,0) to [out=240, in=300] (2.5, 1.44338);
\draw (2.5, 1.44338) to [out=60,in=0] (2.5,4.3301) to [out=180, in=120] (2.5, 1.44338);
\end{tikzpicture}
\begin{tikzpicture}[scale=0.8]
\draw  (0,0) node (1) {} -- (2.5,0) node (2) {} -- (1.25,2.1650635094611) node (3) {}--cycle;\fill  (0,0) circle[radius=1.3pt]  (2.5,0)  circle[radius=1.3pt]  (1.25,2.1650635094611) circle[radius=1.3pt] ;
\draw (0,0) to  [out=180,in=120]  (-.433,-.25)  to [out=300, in=240] (0,0);
\draw (2.5,0) to [out=0,in=60]   (2.933,-.25)  to [out=240, in=300] (2.5,0);
\draw (1.25,2.16506) to  [out=60,in=0]  (1.25,2.66506)  to  [out=180, in=120] (1.25,2.16506);
\end{tikzpicture}
\begin{tikzpicture}[scale=0.8]
\draw (1,0) node{} -- (1.5,0);\draw (0.5,0.86602540378444) node{} -- (0.75,1.2990381056767);\draw (2,0.86602540378444) node{} -- (1.75,1.2990381056767);\fill  (1,0) circle[radius=1.3pt]  (1.5,0)  circle[radius=1.3pt]  (0.5,0.86602540378444) circle[radius=1.3pt]  (0.75,1.2990381056767) circle[radius=1.3pt]   (2,0.86602540378444) circle[radius=1.3pt]  (1.75,1.2990381056767) circle[radius=1.3pt]     ;\draw  (0,0) node (11) {} -- (1,0) node (12) {} -- (0.5,0.86602540378444) node (13) {}--cycle;\fill  (0,0) circle[radius=1.3pt]  (1,0)  circle[radius=1.3pt]  (0.5,0.86602540378444) circle[radius=1.3pt] ;\draw  (1.5,0) node (21) {} -- (2.5,0) node (22) {} -- (2,0.86602540378444) node (23) {}--cycle;\fill  (1.5,0) circle[radius=1.3pt]  (2.5,0)  circle[radius=1.3pt]  (2,0.86602540378444) circle[radius=1.3pt] ;\draw  (0.75,1.2990381056767) node (31) {} -- (1.75,1.2990381056767) node (32) {} -- (1.25,2.1650635094611) node (33) {}--cycle;\fill  (0.75,1.2990381056767) circle[radius=1.3pt]  (1.75,1.2990381056767)  circle[radius=1.3pt]  (1.25,2.1650635094611) circle[radius=1.3pt] ;
\draw (0,0) to  [out=180,in=120]  (-.433,-.25)  to [out=300, in=240] (0,0);
\draw (2.5,0) to [out=0,in=60]   (2.933,-.25)  to [out=240, in=300] (2.5,0);
\draw (1.25,2.16506) to  [out=60,in=0]  (1.25,2.66506)  to  [out=180, in=120] (1.25,2.16506);
\end{tikzpicture}
\begin{tikzpicture}[scale=0.8]
\draw (1,0) node{} -- (1.5,0);\draw (0.5,0.86602540378444) node{} -- (0.75,1.2990381056767);\draw (2,0.86602540378444) node{} -- (1.75,1.2990381056767);\fill  (1,0) circle[radius=1.3pt]  (1.5,0)  circle[radius=1.3pt]  (0.5,0.86602540378444) circle[radius=1.3pt]  (0.75,1.2990381056767) circle[radius=1.3pt]   (2,0.86602540378444) circle[radius=1.3pt]  (1.75,1.2990381056767) circle[radius=1.3pt]     ;\draw (0.4,0) node{} -- (0.6,0);\draw (0.2,0.34641016151378) node{} -- (0.3,0.51961524227066);\draw (0.8,0.34641016151378) node{} -- (0.7,0.51961524227066);\fill  (0.4,0) circle[radius=1.3pt]  (0.6,0)  circle[radius=1.3pt]  (0.2,0.34641016151378) circle[radius=1.3pt]  (0.3,0.51961524227066) circle[radius=1.3pt]   (0.8,0.34641016151378) circle[radius=1.3pt]  (0.7,0.51961524227066) circle[radius=1.3pt]     ;\draw  (0,0) node (111) {} -- (0.4,0) node (112) {} -- (0.2,0.34641016151378) node (113) {}--cycle;\fill  (0,0) circle[radius=1.3pt]  (0.4,0)  circle[radius=1.3pt]  (0.2,0.34641016151378) circle[radius=1.3pt] ;\draw  (0.6,0) node (121) {} -- (1,0) node (122) {} -- (0.8,0.34641016151378) node (123) {}--cycle;\fill  (0.6,0) circle[radius=1.3pt]  (1,0)  circle[radius=1.3pt]  (0.8,0.34641016151378) circle[radius=1.3pt] ;\draw  (0.3,0.51961524227066) node (131) {} -- (0.7,0.51961524227066) node (132) {} -- (0.5,0.86602540378444) node (133) {}--cycle;\fill  (0.3,0.51961524227066) circle[radius=1.3pt]  (0.7,0.51961524227066)  circle[radius=1.3pt]  (0.5,0.86602540378444) circle[radius=1.3pt] ;\draw (1.9,0) node{} -- (2.1,0);\draw (1.7,0.34641016151378) node{} -- (1.8,0.51961524227066);\draw (2.3,0.34641016151378) node{} -- (2.2,0.51961524227066);\fill  (1.9,0) circle[radius=1.3pt]  (2.1,0)  circle[radius=1.3pt]  (1.7,0.34641016151378) circle[radius=1.3pt]  (1.8,0.51961524227066) circle[radius=1.3pt]   (2.3,0.34641016151378) circle[radius=1.3pt]  (2.2,0.51961524227066) circle[radius=1.3pt]     ;\draw  (1.5,0) node (211) {} -- (1.9,0) node (212) {} -- (1.7,0.34641016151378) node (213) {}--cycle;\fill  (1.5,0) circle[radius=1.3pt]  (1.9,0)  circle[radius=1.3pt]  (1.7,0.34641016151378) circle[radius=1.3pt] ;\draw  (2.1,0) node (221) {} -- (2.5,0) node (222) {} -- (2.3,0.34641016151378) node (223) {}--cycle;\fill  (2.1,0) circle[radius=1.3pt]  (2.5,0)  circle[radius=1.3pt]  (2.3,0.34641016151378) circle[radius=1.3pt] ;\draw  (1.8,0.51961524227066) node (231) {} -- (2.2,0.51961524227066) node (232) {} -- (2,0.86602540378444) node (233) {}--cycle;\fill  (1.8,0.51961524227066) circle[radius=1.3pt]  (2.2,0.51961524227066)  circle[radius=1.3pt]  (2,0.86602540378444) circle[radius=1.3pt] ;\draw (1.15,1.2990381056767) node{} -- (1.35,1.2990381056767);\draw (0.95,1.6454482671904) node{} -- (1.05,1.8186533479473);\draw (1.55,1.6454482671904) node{} -- (1.45,1.8186533479473);\fill  (1.15,1.2990381056767) circle[radius=1.3pt]  (1.35,1.2990381056767)  circle[radius=1.3pt]  (0.95,1.6454482671904) circle[radius=1.3pt]  (1.05,1.8186533479473) circle[radius=1.3pt]   (1.55,1.6454482671904) circle[radius=1.3pt]  (1.45,1.8186533479473) circle[radius=1.3pt]     ;\draw  (0.75,1.2990381056767) node (311) {} -- (1.15,1.2990381056767) node (312) {} -- (0.95,1.6454482671904) node (313) {}--cycle;\fill  (0.75,1.2990381056767) circle[radius=1.3pt]  (1.15,1.2990381056767)  circle[radius=1.3pt]  (0.95,1.6454482671904) circle[radius=1.3pt] ;\draw  (1.35,1.2990381056767) node (321) {} -- (1.75,1.2990381056767) node (322) {} -- (1.55,1.6454482671904) node (323) {}--cycle;\fill  (1.35,1.2990381056767) circle[radius=1.3pt]  (1.75,1.2990381056767)  circle[radius=1.3pt]  (1.55,1.6454482671904) circle[radius=1.3pt] ;\draw  (1.05,1.8186533479473) node (331) {} -- (1.45,1.8186533479473) node (332) {} -- (1.25,2.1650635094611) node (333) {}--cycle;\fill  (1.05,1.8186533479473) circle[radius=1.3pt]  (1.45,1.8186533479473)  circle[radius=1.3pt]  (1.25,2.1650635094611) circle[radius=1.3pt] ;
\draw (0,0) to  [out=180,in=120]  (-.433,-.25)  to [out=300, in=240] (0,0);
\draw (2.5,0) to [out=0,in=60]   (2.933,-.25)  to [out=240, in=300] (2.5,0);
\draw (1.25,2.16506) to  [out=60,in=0]  (1.25,2.66506)  to  [out=180, in=120] (1.25,2.16506);
\end{tikzpicture}
\caption{The graphs $H_n$, $n=0,1,2,3$.}\label{fig:Hn}
\end{figure}

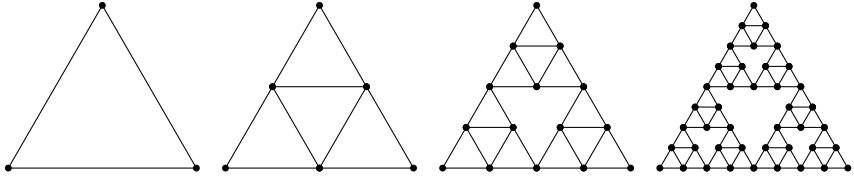
\begin{figure}
\begin{tikzpicture}
\draw  (0,0) node (1) {} -- (2.5,0) node (2) {} -- (1.25,2.1650635094611) node (3) {}--cycle;\fill  (0,0) circle[radius=1.3pt]  (2.5,0)  circle[radius=1.3pt]  (1.25,2.1650635094611) circle[radius=1.3pt] ;
\end{tikzpicture}
\begin{tikzpicture}
\draw  (0,0) node (11) {} -- (1.25,0) node (12) {} -- (0.625,1.0825317547305) node (13) {}--cycle;\fill  (0,0) circle[radius=1.3pt]  (1.25,0)  circle[radius=1.3pt]  (0.625,1.0825317547305) circle[radius=1.3pt] ;\draw  (1.25,0) node (21) {} -- (2.5,0) node (22) {} -- (1.875,1.0825317547305) node (23) {}--cycle;\fill  (1.25,0) circle[radius=1.3pt]  (2.5,0)  circle[radius=1.3pt]  (1.875,1.0825317547305) circle[radius=1.3pt] ;\draw  (0.625,1.0825317547305) node (31) {} -- (1.875,1.0825317547305) node (32) {} -- (1.25,2.1650635094611) node (33) {}--cycle;\fill  (0.625,1.0825317547305) circle[radius=1.3pt]  (1.875,1.0825317547305)  circle[radius=1.3pt]  (1.25,2.1650635094611) circle[radius=1.3pt] ;
\end{tikzpicture}
\begin{tikzpicture}
\draw  (0,0) node (111) {} -- (0.625,0) node (112) {} -- (0.3125,0.54126587736527) node (113) {}--cycle;\fill  (0,0) circle[radius=1.3pt]  (0.625,0)  circle[radius=1.3pt]  (0.3125,0.54126587736527) circle[radius=1.3pt] ;\draw  (0.625,0) node (121) {} -- (1.25,0) node (122) {} -- (0.9375,0.54126587736527) node (123) {}--cycle;\fill  (0.625,0) circle[radius=1.3pt]  (1.25,0)  circle[radius=1.3pt]  (0.9375,0.54126587736527) circle[radius=1.3pt] ;\draw  (0.3125,0.54126587736527) node (131) {} -- (0.9375,0.54126587736527) node (132) {} -- (0.625,1.0825317547305) node (133) {}--cycle;\fill  (0.3125,0.54126587736527) circle[radius=1.3pt]  (0.9375,0.54126587736527)  circle[radius=1.3pt]  (0.625,1.0825317547305) circle[radius=1.3pt] ;\draw  (1.25,0) node (211) {} -- (1.875,0) node (212) {} -- (1.5625,0.54126587736527) node (213) {}--cycle;\fill  (1.25,0) circle[radius=1.3pt]  (1.875,0)  circle[radius=1.3pt]  (1.5625,0.54126587736527) circle[radius=1.3pt] ;\draw  (1.875,0) node (221) {} -- (2.5,0) node (222) {} -- (2.1875,0.54126587736527) node (223) {}--cycle;\fill  (1.875,0) circle[radius=1.3pt]  (2.5,0)  circle[radius=1.3pt]  (2.1875,0.54126587736527) circle[radius=1.3pt] ;\draw  (1.5625,0.54126587736527) node (231) {} -- (2.1875,0.54126587736527) node (232) {} -- (1.875,1.0825317547305) node (233) {}--cycle;\fill  (1.5625,0.54126587736527) circle[radius=1.3pt]  (2.1875,0.54126587736527)  circle[radius=1.3pt]  (1.875,1.0825317547305) circle[radius=1.3pt] ;\draw  (0.625,1.0825317547305) node (311) {} -- (1.25,1.0825317547305) node (312) {} -- (0.9375,1.6237976320958) node (313) {}--cycle;\fill  (0.625,1.0825317547305) circle[radius=1.3pt]  (1.25,1.0825317547305)  circle[radius=1.3pt]  (0.9375,1.6237976320958) circle[radius=1.3pt] ;\draw  (1.25,1.0825317547305) node (321) {} -- (1.875,1.0825317547305) node (322) {} -- (1.5625,1.6237976320958) node (323) {}--cycle;\fill  (1.25,1.0825317547305) circle[radius=1.3pt]  (1.875,1.0825317547305)  circle[radius=1.3pt]  (1.5625,1.6237976320958) circle[radius=1.3pt] ;\draw  (0.9375,1.6237976320958) node (331) {} -- (1.5625,1.6237976320958) node (332) {} -- (1.25,2.1650635094611) node (333) {}--cycle;\fill  (0.9375,1.6237976320958) circle[radius=1.3pt]  (1.5625,1.6237976320958)  circle[radius=1.3pt]  (1.25,2.1650635094611) circle[radius=1.3pt] ;
\end{tikzpicture}
\begin{tikzpicture}
\draw  (0,0) node (1111) {} -- (0.3125,0) node (1112) {} -- (0.15625,0.27063293868264) node (1113) {}--cycle;\fill  (0,0) circle[radius=1.3pt]  (0.3125,0)  circle[radius=1.3pt]  (0.15625,0.27063293868264) circle[radius=1.3pt] ;\draw  (0.3125,0) node (1121) {} -- (0.625,0) node (1122) {} -- (0.46875,0.27063293868264) node (1123) {}--cycle;\fill  (0.3125,0) circle[radius=1.3pt]  (0.625,0)  circle[radius=1.3pt]  (0.46875,0.27063293868264) circle[radius=1.3pt] ;\draw  (0.15625,0.27063293868264) node (1131) {} -- (0.46875,0.27063293868264) node (1132) {} -- (0.3125,0.54126587736527) node (1133) {}--cycle;\fill  (0.15625,0.27063293868264) circle[radius=1.3pt]  (0.46875,0.27063293868264)  circle[radius=1.3pt]  (0.3125,0.54126587736527) circle[radius=1.3pt] ;\draw  (0.625,0) node (1211) {} -- (0.9375,0) node (1212) {} -- (0.78125,0.27063293868264) node (1213) {}--cycle;\fill  (0.625,0) circle[radius=1.3pt]  (0.9375,0)  circle[radius=1.3pt]  (0.78125,0.27063293868264) circle[radius=1.3pt] ;\draw  (0.9375,0) node (1221) {} -- (1.25,0) node (1222) {} -- (1.09375,0.27063293868264) node (1223) {}--cycle;\fill  (0.9375,0) circle[radius=1.3pt]  (1.25,0)  circle[radius=1.3pt]  (1.09375,0.27063293868264) circle[radius=1.3pt] ;\draw  (0.78125,0.27063293868264) node (1231) {} -- (1.09375,0.27063293868264) node (1232) {} -- (0.9375,0.54126587736527) node (1233) {}--cycle;\fill  (0.78125,0.27063293868264) circle[radius=1.3pt]  (1.09375,0.27063293868264)  circle[radius=1.3pt]  (0.9375,0.54126587736527) circle[radius=1.3pt] ;\draw  (0.3125,0.54126587736527) node (1311) {} -- (0.625,0.54126587736527) node (1312) {} -- (0.46875,0.81189881604791) node (1313) {}--cycle;\fill  (0.3125,0.54126587736527) circle[radius=1.3pt]  (0.625,0.54126587736527)  circle[radius=1.3pt]  (0.46875,0.81189881604791) circle[radius=1.3pt] ;\draw  (0.625,0.54126587736527) node (1321) {} -- (0.9375,0.54126587736527) node (1322) {} -- (0.78125,0.81189881604791) node (1323) {}--cycle;\fill  (0.625,0.54126587736527) circle[radius=1.3pt]  (0.9375,0.54126587736527)  circle[radius=1.3pt]  (0.78125,0.81189881604791) circle[radius=1.3pt] ;\draw  (0.46875,0.81189881604791) node (1331) {} -- (0.78125,0.81189881604791) node (1332) {} -- (0.625,1.0825317547305) node (1333) {}--cycle;\fill  (0.46875,0.81189881604791) circle[radius=1.3pt]  (0.78125,0.81189881604791)  circle[radius=1.3pt]  (0.625,1.0825317547305) circle[radius=1.3pt] ;\draw  (1.25,0) node (2111) {} -- (1.5625,0) node (2112) {} -- (1.40625,0.27063293868264) node (2113) {}--cycle;\fill  (1.25,0) circle[radius=1.3pt]  (1.5625,0)  circle[radius=1.3pt]  (1.40625,0.27063293868264) circle[radius=1.3pt] ;\draw  (1.5625,0) node (2121) {} -- (1.875,0) node (2122) {} -- (1.71875,0.27063293868264) node (2123) {}--cycle;\fill  (1.5625,0) circle[radius=1.3pt]  (1.875,0)  circle[radius=1.3pt]  (1.71875,0.27063293868264) circle[radius=1.3pt] ;\draw  (1.40625,0.27063293868264) node (2131) {} -- (1.71875,0.27063293868264) node (2132) {} -- (1.5625,0.54126587736527) node (2133) {}--cycle;\fill  (1.40625,0.27063293868264) circle[radius=1.3pt]  (1.71875,0.27063293868264)  circle[radius=1.3pt]  (1.5625,0.54126587736527) circle[radius=1.3pt] ;\draw  (1.875,0) node (2211) {} -- (2.1875,0) node (2212) {} -- (2.03125,0.27063293868264) node (2213) {}--cycle;\fill  (1.875,0) circle[radius=1.3pt]  (2.1875,0)  circle[radius=1.3pt]  (2.03125,0.27063293868264) circle[radius=1.3pt] ;\draw  (2.1875,0) node (2221) {} -- (2.5,0) node (2222) {} -- (2.34375,0.27063293868264) node (2223) {}--cycle;\fill  (2.1875,0) circle[radius=1.3pt]  (2.5,0)  circle[radius=1.3pt]  (2.34375,0.27063293868264) circle[radius=1.3pt] ;\draw  (2.03125,0.27063293868264) node (2231) {} -- (2.34375,0.27063293868264) node (2232) {} -- (2.1875,0.54126587736527) node (2233) {}--cycle;\fill  (2.03125,0.27063293868264) circle[radius=1.3pt]  (2.34375,0.27063293868264)  circle[radius=1.3pt]  (2.1875,0.54126587736527) circle[radius=1.3pt] ;\draw  (1.5625,0.54126587736527) node (2311) {} -- (1.875,0.54126587736527) node (2312) {} -- (1.71875,0.81189881604791) node (2313) {}--cycle;\fill  (1.5625,0.54126587736527) circle[radius=1.3pt]  (1.875,0.54126587736527)  circle[radius=1.3pt]  (1.71875,0.81189881604791) circle[radius=1.3pt] ;\draw  (1.875,0.54126587736527) node (2321) {} -- (2.1875,0.54126587736527) node (2322) {} -- (2.03125,0.81189881604791) node (2323) {}--cycle;\fill  (1.875,0.54126587736527) circle[radius=1.3pt]  (2.1875,0.54126587736527)  circle[radius=1.3pt]  (2.03125,0.81189881604791) circle[radius=1.3pt] ;\draw  (1.71875,0.81189881604791) node (2331) {} -- (2.03125,0.81189881604791) node (2332) {} -- (1.875,1.0825317547305) node (2333) {}--cycle;\fill  (1.71875,0.81189881604791) circle[radius=1.3pt]  (2.03125,0.81189881604791)  circle[radius=1.3pt]  (1.875,1.0825317547305) circle[radius=1.3pt] ;\draw  (0.625,1.0825317547305) node (3111) {} -- (0.9375,1.0825317547305) node (3112) {} -- (0.78125,1.3531646934132) node (3113) {}--cycle;\fill  (0.625,1.0825317547305) circle[radius=1.3pt]  (0.9375,1.0825317547305)  circle[radius=1.3pt]  (0.78125,1.3531646934132) circle[radius=1.3pt] ;\draw  (0.9375,1.0825317547305) node (3121) {} -- (1.25,1.0825317547305) node (3122) {} -- (1.09375,1.3531646934132) node (3123) {}--cycle;\fill  (0.9375,1.0825317547305) circle[radius=1.3pt]  (1.25,1.0825317547305)  circle[radius=1.3pt]  (1.09375,1.3531646934132) circle[radius=1.3pt] ;\draw  (0.78125,1.3531646934132) node (3131) {} -- (1.09375,1.3531646934132) node (3132) {} -- (0.9375,1.6237976320958) node (3133) {}--cycle;\fill  (0.78125,1.3531646934132) circle[radius=1.3pt]  (1.09375,1.3531646934132)  circle[radius=1.3pt]  (0.9375,1.6237976320958) circle[radius=1.3pt] ;\draw  (1.25,1.0825317547305) node (3211) {} -- (1.5625,1.0825317547305) node (3212) {} -- (1.40625,1.3531646934132) node (3213) {}--cycle;\fill  (1.25,1.0825317547305) circle[radius=1.3pt]  (1.5625,1.0825317547305)  circle[radius=1.3pt]  (1.40625,1.3531646934132) circle[radius=1.3pt] ;\draw  (1.5625,1.0825317547305) node (3221) {} -- (1.875,1.0825317547305) node (3222) {} -- (1.71875,1.3531646934132) node (3223) {}--cycle;\fill  (1.5625,1.0825317547305) circle[radius=1.3pt]  (1.875,1.0825317547305)  circle[radius=1.3pt]  (1.71875,1.3531646934132) circle[radius=1.3pt] ;\draw  (1.40625,1.3531646934132) node (3231) {} -- (1.71875,1.3531646934132) node (3232) {} -- (1.5625,1.6237976320958) node (3233) {}--cycle;\fill  (1.40625,1.3531646934132) circle[radius=1.3pt]  (1.71875,1.3531646934132)  circle[radius=1.3pt]  (1.5625,1.6237976320958) circle[radius=1.3pt] ;\draw  (0.9375,1.6237976320958) node (3311) {} -- (1.25,1.6237976320958) node (3312) {} -- (1.09375,1.8944305707785) node (3313) {}--cycle;\fill  (0.9375,1.6237976320958) circle[radius=1.3pt]  (1.25,1.6237976320958)  circle[radius=1.3pt]  (1.09375,1.8944305707785) circle[radius=1.3pt] ;\draw  (1.25,1.6237976320958) node (3321) {} -- (1.5625,1.6237976320958) node (3322) {} -- (1.40625,1.8944305707785) node (3323) {}--cycle;\fill  (1.25,1.6237976320958) circle[radius=1.3pt]  (1.5625,1.6237976320958)  circle[radius=1.3pt]  (1.40625,1.8944305707785) circle[radius=1.3pt] ;\draw  (1.09375,1.8944305707785) node (3331) {} -- (1.40625,1.8944305707785) node (3332) {} -- (1.25,2.1650635094611) node (3333) {}--cycle;\fill  (1.09375,1.8944305707785) circle[radius=1.3pt]  (1.40625,1.8944305707785)  circle[radius=1.3pt]  (1.25,2.1650635094611) circle[radius=1.3pt] ;
\end{tikzpicture}

\caption{The graphs $G_n$, $n=0,1,2,3$.}\label{fig:Gn}
\end{figure}

We consider the spectra of the Laplacians for the sequence of graphs $H_n$ in Figure~\ref{fig:Hn}, which are the Schreier graphs of the action of the Hanoi towers group on three pegs on the $3$-regular rooted tree, and the spectra of the Laplacians for the sequence of Sierpinski gasket graphs $G_n$ in Figure~\ref{fig:Gn}.  Both sequences limit to the Sierpinski gasket in a natural manner, see~\cite{StrichartzAverages,Nekbook} for $H_n$ and~\cite{Kigbook} for $G_n$.

On the graphs $G_n$ the spectrum of the Laplacian $\Delta^G_n$ may be obtained from that of $\Delta^G_{n-1}$ by a method called spectral decimation~\cite{RammalToulouse,FukushimaShima,MT}; the decimation map involved is independent of $n$, so iteration of it provides a dynamical system for computing $\sigma(\Delta^G_n)$, see Theorem~\ref{thm:SGspect} below.  This technique also gives access to spectra of structures based on these graphs~\cite{StrichartzFractafolds,StrichartzTeplyaev}.   An independently discovered approach applies to the graphs $H_n$ arising from their realization as the Schreier graphs of the self-similar Hanoi towers group, see~\cite{GrigorchukSunic}. In this latter setting one first obtains a two-dimensional dynamics for the spectrum and then proves it is semi-conjugate to a one-dimensional dynamical system; again the spectrum $\sigma(\Delta^H_n)$ can be recovered from the resulting simple dynamics.

In~\cite{GrigorchukSunic} it is suggested that their main result might be derived by spectral decimation, and indeed it is easy to check that the dynamical systems for $G_n$ and $H_n$ are the same. Our purpose here is to examine the connection between the graphs $G_n$ and $H_n$ in more detail, use this to explain a connection between the methods used to obtain the dynamics for the spectra, and determine the structure of eigenfunctions on the graphs $H_n$.   The main result is as follows; definitions of the Laplacian operators are in Section~\ref{sec:graphandLapl}.

\begin{theorem}\label{thm:estructure of Hn}
The spectrum of $\Delta^H_n$ is the set
\begin{equation*}
	\sigma(\Delta^H_n) = \{0\} \cup\biggl( \bigcup_{i=0}^{n-1} R_H^{-i}(-1) \biggr) \cup \biggl( \bigcup_{j=0}^{n-2} R_H^{-j}\Bigl(-\frac53 \Bigr)\biggr)
\end{equation*}
where $R_H(z)=z(3z+5)$. 
The multiplicities and eigenspaces are as follows:
\begin{itemize}
\item $0$ has multiplicity $1$ with constant eigenfunction
\item  $R_H^{-i}(-1)$ contains $2^i$ eigenvalues when $0\leq i\leq n-1$, each with multiplicity $\frac12(3^{n-1-i}+3)$.  For each eigenvalue there is a basis for the eigenspace in which basis elements are localized to the union of two adjacent copies of $H_{i+2}$ in $H_n$, including copies that are self-adjacent via a loop.  (In the case $i=n-1$ this description is vacuous and the support is all of $H_n$.) The basis elements can be derived from a specific set of values on $H_{n-i}$, see Figure~\ref{fig:Hn-1efn}, by $i$ iterations of an explicit algorithm~\eqref{eq:Hnextension}.
\item $R_H^{-j}(-\frac53)$ contains $2^j$ eigenvalues when $0\leq j\leq n-2$, each with multiplicity $\frac12(3^{n-1-j}-1)$. Eigenfunctions reflect the homology of $H_{n-1-j}$: a basis may be obtained by taking alternating values around the corresponding cycles in $H_{n-j}$ and extending to $H_k$ for $n-j+1\leq k\leq n$ by iteration of an explicit algorithm~\eqref{eq:Hnextension}.
\end{itemize}
\end{theorem}

It should be noted that the description of the spectrum is equivalent to that given in Theorem~1.1 of~\cite{GrigorchukSunic}, but the description of the eigenfunctions is new, as is our method. In particular we believe the topological content of the spectrum has not previously been described, though it is substantially the same as for the Sierpinski gasket graphs.  A comparison of our method to that of~\cite{GrigorchukSunic} is in Section~\ref{sec:GScompare}.  The latter has also been used to treat other important examples of Schreier graphs from self-similar groups~\cite{BG}.

The structure of the paper is as follows. In Section~\ref{sec:graphandLapl} we give formal definitions of our graphs and Laplacians. Section~\ref{sec:specdec}  contains the definition and some known results about spectral similarity, which is our main tool.  Our results are in Section~\ref{sec:specsims}, which establishes spectral similarities between the Laplacians on the graphs $G_n$ and $H_n$, Section~\ref{sec:specselfsim}, where we prove  the Laplacian on $H_n$ is spectrally self-similar and establish the iteration~\eqref{eq:Hnextension} for obtaining eigenfunctions of $H_{n+1}$ from those of $H_n$, and Section~\ref{sec:efns}, where we discuss the construction of eigenfunctions and prove Theorem~\ref{thm:estructure of Hn}.


\section{Graphs and Laplacians}\label{sec:graphandLapl}

\subsection*{Graphs}

A rapid way to formally define the graphs $G_n$ in Figure~\ref{fig:Gn} is to take $V^G_0=\{p_1,p_2,p_3\}$ to be the vertices of a triangle in the plane and let $G_0$ be the complete graph on $V^G_0$. Then set  $F_j(x)=\frac12(x+p_j)$ for $j=1,2,3$ and $V^G_n=\cup_{j=1}^3F_j(V_{n-1})$ for $n\geq1$. $G_n$ has vertex set $V^G_n$ and inductively defined edge relation $x\sim_{G_n}y$ if and only if $x=F_j(x')$, $y=F_j(y')$ for  some $j\in\{1,2,3\}$ and points $x'\sim_{G_{n-1}}y'$ from $V^G_{n-1}$. Figure~\ref{fig:Gn} is then derived by iterative application of the three maps $F_j$ to the graph $G_0$. Note that the three copies $F_j(V^G_{n-1})$ overlap at three points, so there are $\frac12(3^n+3)$ points in $V^G_n$.

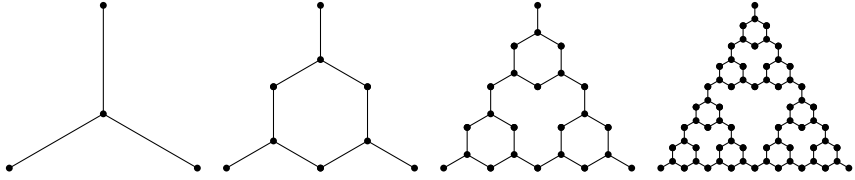
\begin{figure}
\begin{tikzpicture}
\draw  (0,0) node (1) {} -- (1.25,0.72168783648703) node (0) {} -- (2.5,0) node (2) {}; \draw  (1.25,0.72168783648703)  -- (1.25,2.1650635094611) node (3) {}; \fill  (0,0) circle[radius=1.3pt]  (2.5,0)  circle[radius=1.3pt]  (1.25,2.1650635094611) circle[radius=1.3pt]  (1.25,0.72168783648703) circle[radius=1.3pt]  ;
\end{tikzpicture}
\begin{tikzpicture}
\draw  (0,0) node (11) {} -- (0.625,0.36084391824352) node (10) {} -- (1.25,0) node (12) {}; \draw  (0.625,0.36084391824352)  -- (0.625,1.0825317547305) node (13) {}; \fill  (0,0) circle[radius=1.3pt]  (1.25,0)  circle[radius=1.3pt]  (0.625,1.0825317547305) circle[radius=1.3pt]  (0.625,0.36084391824352) circle[radius=1.3pt]  ;\draw  (1.25,0) node (21) {} -- (1.875,0.36084391824352) node (20) {} -- (2.5,0) node (22) {}; \draw  (1.875,0.36084391824352)  -- (1.875,1.0825317547305) node (23) {}; \fill  (1.25,0) circle[radius=1.3pt]  (2.5,0)  circle[radius=1.3pt]  (1.875,1.0825317547305) circle[radius=1.3pt]  (1.875,0.36084391824352) circle[radius=1.3pt]  ;\draw  (0.625,1.0825317547305) node (31) {} -- (1.25,1.4433756729741) node (30) {} -- (1.875,1.0825317547305) node (32) {}; \draw  (1.25,1.4433756729741)  -- (1.25,2.1650635094611) node (33) {}; \fill  (0.625,1.0825317547305) circle[radius=1.3pt]  (1.875,1.0825317547305)  circle[radius=1.3pt]  (1.25,2.1650635094611) circle[radius=1.3pt]  (1.25,1.4433756729741) circle[radius=1.3pt]  ;
\end{tikzpicture}
\begin{tikzpicture}
\draw  (0,0) node (111) {} -- (0.3125,0.18042195912176) node (110) {} -- (0.625,0) node (112) {}; \draw  (0.3125,0.18042195912176)  -- (0.3125,0.54126587736527) node (113) {}; \fill  (0,0) circle[radius=1.3pt]  (0.625,0)  circle[radius=1.3pt]  (0.3125,0.54126587736527) circle[radius=1.3pt]  (0.3125,0.18042195912176) circle[radius=1.3pt]  ;\draw  (0.625,0) node (121) {} -- (0.9375,0.18042195912176) node (120) {} -- (1.25,0) node (122) {}; \draw  (0.9375,0.18042195912176)  -- (0.9375,0.54126587736527) node (123) {}; \fill  (0.625,0) circle[radius=1.3pt]  (1.25,0)  circle[radius=1.3pt]  (0.9375,0.54126587736527) circle[radius=1.3pt]  (0.9375,0.18042195912176) circle[radius=1.3pt]  ;\draw  (0.3125,0.54126587736527) node (131) {} -- (0.625,0.72168783648703) node (130) {} -- (0.9375,0.54126587736527) node (132) {}; \draw  (0.625,0.72168783648703)  -- (0.625,1.0825317547305) node (133) {}; \fill  (0.3125,0.54126587736527) circle[radius=1.3pt]  (0.9375,0.54126587736527)  circle[radius=1.3pt]  (0.625,1.0825317547305) circle[radius=1.3pt]  (0.625,0.72168783648703) circle[radius=1.3pt]  ;\draw  (1.25,0) node (211) {} -- (1.5625,0.18042195912176) node (210) {} -- (1.875,0) node (212) {}; \draw  (1.5625,0.18042195912176)  -- (1.5625,0.54126587736527) node (213) {}; \fill  (1.25,0) circle[radius=1.3pt]  (1.875,0)  circle[radius=1.3pt]  (1.5625,0.54126587736527) circle[radius=1.3pt]  (1.5625,0.18042195912176) circle[radius=1.3pt]  ;\draw  (1.875,0) node (221) {} -- (2.1875,0.18042195912176) node (220) {} -- (2.5,0) node (222) {}; \draw  (2.1875,0.18042195912176)  -- (2.1875,0.54126587736527) node (223) {}; \fill  (1.875,0) circle[radius=1.3pt]  (2.5,0)  circle[radius=1.3pt]  (2.1875,0.54126587736527) circle[radius=1.3pt]  (2.1875,0.18042195912176) circle[radius=1.3pt]  ;\draw  (1.5625,0.54126587736527) node (231) {} -- (1.875,0.72168783648703) node (230) {} -- (2.1875,0.54126587736527) node (232) {}; \draw  (1.875,0.72168783648703)  -- (1.875,1.0825317547305) node (233) {}; \fill  (1.5625,0.54126587736527) circle[radius=1.3pt]  (2.1875,0.54126587736527)  circle[radius=1.3pt]  (1.875,1.0825317547305) circle[radius=1.3pt]  (1.875,0.72168783648703) circle[radius=1.3pt]  ;\draw  (0.625,1.0825317547305) node (311) {} -- (0.9375,1.2629537138523) node (310) {} -- (1.25,1.0825317547305) node (312) {}; \draw  (0.9375,1.2629537138523)  -- (0.9375,1.6237976320958) node (313) {}; \fill  (0.625,1.0825317547305) circle[radius=1.3pt]  (1.25,1.0825317547305)  circle[radius=1.3pt]  (0.9375,1.6237976320958) circle[radius=1.3pt]  (0.9375,1.2629537138523) circle[radius=1.3pt]  ;\draw  (1.25,1.0825317547305) node (321) {} -- (1.5625,1.2629537138523) node (320) {} -- (1.875,1.0825317547305) node (322) {}; \draw  (1.5625,1.2629537138523)  -- (1.5625,1.6237976320958) node (323) {}; \fill  (1.25,1.0825317547305) circle[radius=1.3pt]  (1.875,1.0825317547305)  circle[radius=1.3pt]  (1.5625,1.6237976320958) circle[radius=1.3pt]  (1.5625,1.2629537138523) circle[radius=1.3pt]  ;\draw  (0.9375,1.6237976320958) node (331) {} -- (1.25,1.8042195912176) node (330) {} -- (1.5625,1.6237976320958) node (332) {}; \draw  (1.25,1.8042195912176)  -- (1.25,2.1650635094611) node (333) {}; \fill  (0.9375,1.6237976320958) circle[radius=1.3pt]  (1.5625,1.6237976320958)  circle[radius=1.3pt]  (1.25,2.1650635094611) circle[radius=1.3pt]  (1.25,1.8042195912176) circle[radius=1.3pt]  ;
\end{tikzpicture}
\begin{tikzpicture}
\draw  (0,0) node (1111) {} -- (0.15625,0.090210979560879) node (1110) {} -- (0.3125,0) node (1112) {}; \draw  (0.15625,0.090210979560879)  -- (0.15625,0.27063293868264) node (1113) {}; \fill  (0,0) circle[radius=1.3pt]  (0.3125,0)  circle[radius=1.3pt]  (0.15625,0.27063293868264) circle[radius=1.3pt]  (0.15625,0.090210979560879) circle[radius=1.3pt]  ;\draw  (0.3125,0) node (1121) {} -- (0.46875,0.090210979560879) node (1120) {} -- (0.625,0) node (1122) {}; \draw  (0.46875,0.090210979560879)  -- (0.46875,0.27063293868264) node (1123) {}; \fill  (0.3125,0) circle[radius=1.3pt]  (0.625,0)  circle[radius=1.3pt]  (0.46875,0.27063293868264) circle[radius=1.3pt]  (0.46875,0.090210979560879) circle[radius=1.3pt]  ;\draw  (0.15625,0.27063293868264) node (1131) {} -- (0.3125,0.36084391824352) node (1130) {} -- (0.46875,0.27063293868264) node (1132) {}; \draw  (0.3125,0.36084391824352)  -- (0.3125,0.54126587736527) node (1133) {}; \fill  (0.15625,0.27063293868264) circle[radius=1.3pt]  (0.46875,0.27063293868264)  circle[radius=1.3pt]  (0.3125,0.54126587736527) circle[radius=1.3pt]  (0.3125,0.36084391824352) circle[radius=1.3pt]  ;\draw  (0.625,0) node (1211) {} -- (0.78125,0.090210979560879) node (1210) {} -- (0.9375,0) node (1212) {}; \draw  (0.78125,0.090210979560879)  -- (0.78125,0.27063293868264) node (1213) {}; \fill  (0.625,0) circle[radius=1.3pt]  (0.9375,0)  circle[radius=1.3pt]  (0.78125,0.27063293868264) circle[radius=1.3pt]  (0.78125,0.090210979560879) circle[radius=1.3pt]  ;\draw  (0.9375,0) node (1221) {} -- (1.09375,0.090210979560879) node (1220) {} -- (1.25,0) node (1222) {}; \draw  (1.09375,0.090210979560879)  -- (1.09375,0.27063293868264) node (1223) {}; \fill  (0.9375,0) circle[radius=1.3pt]  (1.25,0)  circle[radius=1.3pt]  (1.09375,0.27063293868264) circle[radius=1.3pt]  (1.09375,0.090210979560879) circle[radius=1.3pt]  ;\draw  (0.78125,0.27063293868264) node (1231) {} -- (0.9375,0.36084391824352) node (1230) {} -- (1.09375,0.27063293868264) node (1232) {}; \draw  (0.9375,0.36084391824352)  -- (0.9375,0.54126587736527) node (1233) {}; \fill  (0.78125,0.27063293868264) circle[radius=1.3pt]  (1.09375,0.27063293868264)  circle[radius=1.3pt]  (0.9375,0.54126587736527) circle[radius=1.3pt]  (0.9375,0.36084391824352) circle[radius=1.3pt]  ;\draw  (0.3125,0.54126587736527) node (1311) {} -- (0.46875,0.63147685692615) node (1310) {} -- (0.625,0.54126587736527) node (1312) {}; \draw  (0.46875,0.63147685692615)  -- (0.46875,0.81189881604791) node (1313) {}; \fill  (0.3125,0.54126587736527) circle[radius=1.3pt]  (0.625,0.54126587736527)  circle[radius=1.3pt]  (0.46875,0.81189881604791) circle[radius=1.3pt]  (0.46875,0.63147685692615) circle[radius=1.3pt]  ;\draw  (0.625,0.54126587736527) node (1321) {} -- (0.78125,0.63147685692615) node (1320) {} -- (0.9375,0.54126587736527) node (1322) {}; \draw  (0.78125,0.63147685692615)  -- (0.78125,0.81189881604791) node (1323) {}; \fill  (0.625,0.54126587736527) circle[radius=1.3pt]  (0.9375,0.54126587736527)  circle[radius=1.3pt]  (0.78125,0.81189881604791) circle[radius=1.3pt]  (0.78125,0.63147685692615) circle[radius=1.3pt]  ;\draw  (0.46875,0.81189881604791) node (1331) {} -- (0.625,0.90210979560879) node (1330) {} -- (0.78125,0.81189881604791) node (1332) {}; \draw  (0.625,0.90210979560879)  -- (0.625,1.0825317547305) node (1333) {}; \fill  (0.46875,0.81189881604791) circle[radius=1.3pt]  (0.78125,0.81189881604791)  circle[radius=1.3pt]  (0.625,1.0825317547305) circle[radius=1.3pt]  (0.625,0.90210979560879) circle[radius=1.3pt]  ;\draw  (1.25,0) node (2111) {} -- (1.40625,0.090210979560879) node (2110) {} -- (1.5625,0) node (2112) {}; \draw  (1.40625,0.090210979560879)  -- (1.40625,0.27063293868264) node (2113) {}; \fill  (1.25,0) circle[radius=1.3pt]  (1.5625,0)  circle[radius=1.3pt]  (1.40625,0.27063293868264) circle[radius=1.3pt]  (1.40625,0.090210979560879) circle[radius=1.3pt]  ;\draw  (1.5625,0) node (2121) {} -- (1.71875,0.090210979560879) node (2120) {} -- (1.875,0) node (2122) {}; \draw  (1.71875,0.090210979560879)  -- (1.71875,0.27063293868264) node (2123) {}; \fill  (1.5625,0) circle[radius=1.3pt]  (1.875,0)  circle[radius=1.3pt]  (1.71875,0.27063293868264) circle[radius=1.3pt]  (1.71875,0.090210979560879) circle[radius=1.3pt]  ;\draw  (1.40625,0.27063293868264) node (2131) {} -- (1.5625,0.36084391824352) node (2130) {} -- (1.71875,0.27063293868264) node (2132) {}; \draw  (1.5625,0.36084391824352)  -- (1.5625,0.54126587736527) node (2133) {}; \fill  (1.40625,0.27063293868264) circle[radius=1.3pt]  (1.71875,0.27063293868264)  circle[radius=1.3pt]  (1.5625,0.54126587736527) circle[radius=1.3pt]  (1.5625,0.36084391824352) circle[radius=1.3pt]  ;\draw  (1.875,0) node (2211) {} -- (2.03125,0.090210979560879) node (2210) {} -- (2.1875,0) node (2212) {}; \draw  (2.03125,0.090210979560879)  -- (2.03125,0.27063293868264) node (2213) {}; \fill  (1.875,0) circle[radius=1.3pt]  (2.1875,0)  circle[radius=1.3pt]  (2.03125,0.27063293868264) circle[radius=1.3pt]  (2.03125,0.090210979560879) circle[radius=1.3pt]  ;\draw  (2.1875,0) node (2221) {} -- (2.34375,0.090210979560879) node (2220) {} -- (2.5,0) node (2222) {}; \draw  (2.34375,0.090210979560879)  -- (2.34375,0.27063293868264) node (2223) {}; \fill  (2.1875,0) circle[radius=1.3pt]  (2.5,0)  circle[radius=1.3pt]  (2.34375,0.27063293868264) circle[radius=1.3pt]  (2.34375,0.090210979560879) circle[radius=1.3pt]  ;\draw  (2.03125,0.27063293868264) node (2231) {} -- (2.1875,0.36084391824352) node (2230) {} -- (2.34375,0.27063293868264) node (2232) {}; \draw  (2.1875,0.36084391824352)  -- (2.1875,0.54126587736527) node (2233) {}; \fill  (2.03125,0.27063293868264) circle[radius=1.3pt]  (2.34375,0.27063293868264)  circle[radius=1.3pt]  (2.1875,0.54126587736527) circle[radius=1.3pt]  (2.1875,0.36084391824352) circle[radius=1.3pt]  ;\draw  (1.5625,0.54126587736527) node (2311) {} -- (1.71875,0.63147685692615) node (2310) {} -- (1.875,0.54126587736527) node (2312) {}; \draw  (1.71875,0.63147685692615)  -- (1.71875,0.81189881604791) node (2313) {}; \fill  (1.5625,0.54126587736527) circle[radius=1.3pt]  (1.875,0.54126587736527)  circle[radius=1.3pt]  (1.71875,0.81189881604791) circle[radius=1.3pt]  (1.71875,0.63147685692615) circle[radius=1.3pt]  ;\draw  (1.875,0.54126587736527) node (2321) {} -- (2.03125,0.63147685692615) node (2320) {} -- (2.1875,0.54126587736527) node (2322) {}; \draw  (2.03125,0.63147685692615)  -- (2.03125,0.81189881604791) node (2323) {}; \fill  (1.875,0.54126587736527) circle[radius=1.3pt]  (2.1875,0.54126587736527)  circle[radius=1.3pt]  (2.03125,0.81189881604791) circle[radius=1.3pt]  (2.03125,0.63147685692615) circle[radius=1.3pt]  ;\draw  (1.71875,0.81189881604791) node (2331) {} -- (1.875,0.90210979560879) node (2330) {} -- (2.03125,0.81189881604791) node (2332) {}; \draw  (1.875,0.90210979560879)  -- (1.875,1.0825317547305) node (2333) {}; \fill  (1.71875,0.81189881604791) circle[radius=1.3pt]  (2.03125,0.81189881604791)  circle[radius=1.3pt]  (1.875,1.0825317547305) circle[radius=1.3pt]  (1.875,0.90210979560879) circle[radius=1.3pt]  ;\draw  (0.625,1.0825317547305) node (3111) {} -- (0.78125,1.1727427342914) node (3110) {} -- (0.9375,1.0825317547305) node (3112) {}; \draw  (0.78125,1.1727427342914)  -- (0.78125,1.3531646934132) node (3113) {}; \fill  (0.625,1.0825317547305) circle[radius=1.3pt]  (0.9375,1.0825317547305)  circle[radius=1.3pt]  (0.78125,1.3531646934132) circle[radius=1.3pt]  (0.78125,1.1727427342914) circle[radius=1.3pt]  ;\draw  (0.9375,1.0825317547305) node (3121) {} -- (1.09375,1.1727427342914) node (3120) {} -- (1.25,1.0825317547305) node (3122) {}; \draw  (1.09375,1.1727427342914)  -- (1.09375,1.3531646934132) node (3123) {}; \fill  (0.9375,1.0825317547305) circle[radius=1.3pt]  (1.25,1.0825317547305)  circle[radius=1.3pt]  (1.09375,1.3531646934132) circle[radius=1.3pt]  (1.09375,1.1727427342914) circle[radius=1.3pt]  ;\draw  (0.78125,1.3531646934132) node (3131) {} -- (0.9375,1.4433756729741) node (3130) {} -- (1.09375,1.3531646934132) node (3132) {}; \draw  (0.9375,1.4433756729741)  -- (0.9375,1.6237976320958) node (3133) {}; \fill  (0.78125,1.3531646934132) circle[radius=1.3pt]  (1.09375,1.3531646934132)  circle[radius=1.3pt]  (0.9375,1.6237976320958) circle[radius=1.3pt]  (0.9375,1.4433756729741) circle[radius=1.3pt]  ;\draw  (1.25,1.0825317547305) node (3211) {} -- (1.40625,1.1727427342914) node (3210) {} -- (1.5625,1.0825317547305) node (3212) {}; \draw  (1.40625,1.1727427342914)  -- (1.40625,1.3531646934132) node (3213) {}; \fill  (1.25,1.0825317547305) circle[radius=1.3pt]  (1.5625,1.0825317547305)  circle[radius=1.3pt]  (1.40625,1.3531646934132) circle[radius=1.3pt]  (1.40625,1.1727427342914) circle[radius=1.3pt]  ;\draw  (1.5625,1.0825317547305) node (3221) {} -- (1.71875,1.1727427342914) node (3220) {} -- (1.875,1.0825317547305) node (3222) {}; \draw  (1.71875,1.1727427342914)  -- (1.71875,1.3531646934132) node (3223) {}; \fill  (1.5625,1.0825317547305) circle[radius=1.3pt]  (1.875,1.0825317547305)  circle[radius=1.3pt]  (1.71875,1.3531646934132) circle[radius=1.3pt]  (1.71875,1.1727427342914) circle[radius=1.3pt]  ;\draw  (1.40625,1.3531646934132) node (3231) {} -- (1.5625,1.4433756729741) node (3230) {} -- (1.71875,1.3531646934132) node (3232) {}; \draw  (1.5625,1.4433756729741)  -- (1.5625,1.6237976320958) node (3233) {}; \fill  (1.40625,1.3531646934132) circle[radius=1.3pt]  (1.71875,1.3531646934132)  circle[radius=1.3pt]  (1.5625,1.6237976320958) circle[radius=1.3pt]  (1.5625,1.4433756729741) circle[radius=1.3pt]  ;\draw  (0.9375,1.6237976320958) node (3311) {} -- (1.09375,1.7140086116567) node (3310) {} -- (1.25,1.6237976320958) node (3312) {}; \draw  (1.09375,1.7140086116567)  -- (1.09375,1.8944305707785) node (3313) {}; \fill  (0.9375,1.6237976320958) circle[radius=1.3pt]  (1.25,1.6237976320958)  circle[radius=1.3pt]  (1.09375,1.8944305707785) circle[radius=1.3pt]  (1.09375,1.7140086116567) circle[radius=1.3pt]  ;\draw  (1.25,1.6237976320958) node (3321) {} -- (1.40625,1.7140086116567) node (3320) {} -- (1.5625,1.6237976320958) node (3322) {}; \draw  (1.40625,1.7140086116567)  -- (1.40625,1.8944305707785) node (3323) {}; \fill  (1.25,1.6237976320958) circle[radius=1.3pt]  (1.5625,1.6237976320958)  circle[radius=1.3pt]  (1.40625,1.8944305707785) circle[radius=1.3pt]  (1.40625,1.7140086116567) circle[radius=1.3pt]  ;\draw  (1.09375,1.8944305707785) node (3331) {} -- (1.25,1.9846415503393) node (3330) {} -- (1.40625,1.8944305707785) node (3332) {}; \draw  (1.25,1.9846415503393)  -- (1.25,2.1650635094611) node (3333) {}; \fill  (1.09375,1.8944305707785) circle[radius=1.3pt]  (1.40625,1.8944305707785)  circle[radius=1.3pt]  (1.25,2.1650635094611) circle[radius=1.3pt]  (1.25,1.9846415503393) circle[radius=1.3pt]  ;
\end{tikzpicture}

\caption{The graphs $J_n$, $n=0,1,2,3$.}\label{fig:Jn}
\end{figure}

We define the graphs $H_n$ in Figure~\ref{fig:Hn} indirectly, obtaining along the way a third series of graphs $J_n$, shown in Figure~\ref{fig:Jn} that will be useful in their own right. Let $p_0=\frac13\sum_1^3p_j$ be the midpoint of the previously-defined triangle and  $V^J_0=\{p_0,p_1,p_2,p_3\}$. Define $J_0$ to be the graph on $V^J_0$ with edges from $p_0$ to each $p_j$, $j=1,2,3$.  Then define $J_n$ to have vertex set  $V^J_n=\cup_1^3F_j(V_{n-1})$ for $n\geq1$ and edges $x\sim_{J_n}y$ precisely when $x=F_j(x')$, $y=F_j(y')$ for  some $j\in\{1,2,3\}$ and points $x'\sim_{J_{n-1}}y'$ from $V^J_{n-1}$.

Now we can define $H_n$ on the vertex set obtained inductively from $V^H_0=\{p_0\}$ and $V^H_n=\cup_1^3F_j(V^H_{n-1})$ for $n\geq1$ by providing an edge $x\sim_{H_n}x'$ for $x\neq x'\in V^H_n$ when there is a vertex $y\in V^J_n$ so both $y\sim_{J_n}x$ and $y\sim_{J_n}x'$. Evidently $V^J_n=V^G_n\cup V^H_n$ and this  point $y$ is actually from $V^G_n$.   Moreover for each $y\in V^G_n\setminus V^G_0$ there is a uniquely determined pair $x\neq x'$ with this property. In the case $y\in V^G_0$ we follow~\cite{GrigorchukSunic} and define a special case, placing a loop at the unique neighboring $x\in V^H_n$. The purpose of this loop will become more clear later; for now we note that it will make $H_n$ a $3$-regular graph.

The reader might find this to be a convenient point to check that $H_n$ is the level~$n$ Schreier graph of the action of the Hanoi Towers group on the rooted $3$-tree, as described in~\cite{GrigorchukSunic}. Indeed, if we introduce a word space $W_n=\{1,2,3\}^n$ and for $w=w_1\dotsm w_n\in W_n$ write $F_w=F_{w_1}\circ\dotsm\circ F_{w_n}$ then $w\mapsto F_w(p_0)$ is a bijection $W_n\to V^H_n$ and there is an edge from $F_w(p_0)$ to $F_{w'}(p_0)$ exactly when there are $j,j'\in\{1,2,3\}$ so $F_w(p_j)=F_{w'}(p_{j'})$, which is precisely when the Hanoi Towers group contains an element taking $w$ to $w'$. A characterization of the connection between sequences of self-similar graphs and sequences of Schreier graphs of self-similar groups is in~\cite{KSW}; other connections between the methods used to analyze these sequences are in~\cite{NekTep}.

\subsection*{Laplacian}
There are several standard notions of Laplacian associated to a graph.  In order to simplify the use of some results from~\cite{MT} we use the graph Laplacian, which is defined for a function $f$ on the vertex set by
\begin{equation*}
	\Delta_nf(x)=\frac1{\deg(x)}\sum_{x\sim_n y} (f(y)-f(x))
	\end{equation*}
 where the degree $\deg(x)$ is the number of edges incident at the vertex $x$. However the case of a loop at $x$ is slightly different: then $x=y$ and $f(x)-f(y)=0$, and we count the loop as contributing $1$ to the degree.  It is easy to check that this is consistent with the usual simple random walk in which a mass at $x$ is divided by $\deg(x)$ and propagated to the neighbors, one of which neighbors is $x$ itself if there is a loop at $x$.

Observe that on $G_n$ we have $\deg(x)=4$ for $x\in V^G_n\setminus V^G_0$ and $\deg(x)=2$ on $V^G_0$.  The graph $H_n$ is regular, with $\deg(x)=3$ at all vertices, because the loops contribute $1$ to the degree of the three vertices at the corners of the triangle.  The graphs $J_n$ have $\deg(x)=3$ at those vertices in $V^J_n$ that come from $V^H_n$, $\deg(x)=2$ on $V^G_n\setminus V^G_0$ and $\deg(x)=1$ on $V^G_0$.

Rather than giving a recursion for the entries of the Laplacian matrices, which is possible but would not be used later, we illustrate our definition with the following explicit examples. For clarity we have made blocks of vertices with the same degree and ordered them by decreasing degree. Some other ordering choices have been used to emphasize the symmetry of the matrices.  Setting
\begin{equation*}
	M=\begin{bmatrix}
		0& 1&1 \\
		1 & 0& 1\\
		1 & 1 & 0
		\end{bmatrix}
	\end{equation*}
and writing $I$ for the $2\times3$ identity matrix  and $0$ for the zero matrix we have $\Delta^H_1= \frac13M-I$ and in block form:
\begin{equation*}
	\Delta^G_1= \begin{bmatrix} \frac14 M - I & \frac14 M\\ \frac12M & -I \end{bmatrix} \qquad\qquad
	\Delta^J_1 = \begin{bmatrix} -I & \frac13M & \frac13 M \\ \frac12M & -I & 0 \\ I& 0 & -I \end{bmatrix}.
	\end{equation*} 

\section{Spectral similarity and decimation}\label{sec:specdec}

We present the elementary features of spectral similarity and spectral decimation as they apply to the circumstances in the problem at hand.  These include the connection between the classical Schur complement formula for relating the spectrum of a matrix to that of submatrices in a block decomposition and the notion of spectral similarity, and a construction from~\cite{MT} which permits one to glue spectrally similar graph Laplacians in a manner that preserves spectral similarity.

Let $\mathcal{H}$ be a finite dimensional complex vector space and $\tilde{\mathcal{H}}$ be a subspace. Write $P:\mathcal{H}\to\tilde{\mathcal{H}}$ for the orthogonal projection (so $P^\ast$ is the inclusion) and $Q$ for the projection to the orthogonal complement. Let $D$ be an operator on $\mathcal{H}$ and $\tilde D$ be an operator on $\tilde{\mathcal{H}}$.

In accordance with Definition~2.1 in~\cite{MT} we say $D$ is spectrally similar to $\tilde D$ with functions $\phi$ and $\tilde \phi$ if $P^\ast(D-z)^{-1}P=\bigl( \tilde \phi(z)\tilde D - \phi(z)\bigr)^{-1}$ as meromorphic functions.  This is related to the Schur complement in the following manner (Lemma~3.3 in~\cite{MT}).  Writing the symmetric matrix $D$ in the block form on $\tilde{\mathcal{H}}^{\perp}\otimes\tilde{\mathcal{H}}$
\begin{equation*}
	D = \begin{pmatrix}  QDQ^\ast & PDQ^\ast\\  QDP^\ast & PDP^\ast\\   \end{pmatrix}
	\end{equation*}
we find that $D$ is spectrally similar to $\tilde D$ if and only if the Schur complement from $D-z$ satisfies
\begin{equation*}
	PDP^\ast - z - QDP^\ast (QDQ^\ast - z)^{-1} PDQ^\ast = \tilde\phi(z)\tilde D - \phi(z).
	\end{equation*}
In particular, we can compute the spectrum of $D$ from that of $\tilde D$ using the associated functions.  Writing $\sigma(D)$ for the spectrum we have the following.
\begin{theorem}[\protect{Theorem 3.6 in \cite{MT}}]
\label{thm:spectralcorresp}
If $D$ is spectrally similar to $\tilde{D}$ with functions $\phi$ and $\tilde\phi$, let the exceptional set be $E=\sigma(QDQ^\ast)\cup\{z:\tilde\phi(z)=0\}$ and the (spectral) decimation function be $R(z)=\phi(z)/\tilde\phi(z)$. Suppose $z\notin E$.  Then $R(z)$ is an eigenvalue of $\tilde{D}$ if and only if $z$ is an eigenvalue of $D$, and there is a bijection $\tilde f \mapsto f=\tilde f - (QDQ^\ast -z)^{-1} PDQ^\ast \tilde f$ between the corresponding eigenspaces.
\end{theorem}

Some general theory regarding the exceptional set is summarized in~\cite{BajorinChenJPhs}.

\begin{example}\label{ex:starmap}
A case of special interest in the problem at hand occurs when $D_N$ is the Laplacian of the graph $\Gamma_N$ in which one vertex $p_0$ is connected by a single edge to each of $N$ other vertices $p_j$, $j=1,\dotsc N$, so $\mathcal{H}=\mathbb{C}^{N+1}$. Take $\tilde D_N$ to be the Laplacian of the complete graph on the vertices  $p_j$, $j=1,\dotsc N$ and let $\tilde{\mathcal{H}}$ be the corresponding subspace. Writing the block form 
\begin{equation*}
	D_N-z= \begin{pmatrix}  -(1+z)  & (\frac1N)1_{1\times N}\\  1_{N\times1} &  -(1+z)I_N \end{pmatrix}
	\end{equation*}
where $I_N$ is the $N\times N$ identity and $1_{m\times n}$ is the $m\times n$ matrix with all entries equal $1$ we can compute the Schur complement $S_z$ to be
\begin{equation*}
	S_z = -(1+z)I_N + \frac1{N(1+z)}1_{N\times N} = \frac {N-1}{N(1+z)} \tilde D_N - \Bigl(1+z - \frac1{1+z}\Bigr)I_N
\end{equation*}
so that $D_N$ is spectrally similar to $\tilde D_N$ with $\phi_N(z)=1+z-\frac1{1+z}$ and $\tilde\phi_N(z)=\frac {N-1}{N(1+z)}$. This relation is valid provided $z$ is not in the exceptional set $E$. However, the exceptional set is small: the spectrum of the Schur block is simply $\{-1\}$ and $\tilde\phi_N(z)=\frac {N-1}{N(1+z)}$ is always non-zero, so $E=\{-1\}$ independent of $N$.  The  reader will note that this Schur complement reduction is the spectral version of the star-mesh transform that is well known in electrical network theory; the case $N=3$ of the latter is often called the $\Delta-Y$ transform.

It is worth noting what happens to eigenfunctions in this context. If  $z\neq-1$, $R(z)\in\sigma(\tilde D_N)$ and $\tilde f$ is an eigenfunction of $\tilde D_N$ with eigenvalue $R(z)$, then the corresponding eigenfunction of $D_N$ with eigenvalue $z$ is
\begin{equation*}
	f = \tilde f - (QD_N Q^\ast)^{-1} \tilde f = \tilde f + \frac1{1+z} \bigl(\frac1N\bigr) 1_{1\times N} \tilde f_0
	\end{equation*}
so that $f=f_0$ on $\tilde V$ and $f(p_0) = \frac1{N(1+z)} \sum_{j=1}^N f_0(p_j)$.
The remaining possibility is $z=-1$, in which case the spectral similarity cannot be used but
\begin{equation*}
	D_N+I=\begin{pmatrix} 0 & (\frac1N)1_{1\times N} \\  1_{N\times 1} & 0_{N} \end{pmatrix}
	\end{equation*}
so $-1$ is an eigenvalue with an $N-1$ dimensional eigenspace consisting of those functions $f$ with $f(p_0)=0$ and $\sum_{j=1}^N f(p_j)=0$.
\end{example}

\section{Spectral similarities for $G_n$, $H_n$ and $J_n$}\label{sec:specsims}

The significance of Example	~\ref{ex:starmap} is that the graph $J_n$ can be built by replacing the edges of $H_n$ by copies of $\Gamma_2$, and also by replacing all copies of the graph $V^G_0$ (the complete graph on 3 vertices) in $V^G_n$ with copies of $V^J_0=\Gamma_3$. The effect of this type of graph replacement operation on the spectra of graphs was previously studied in~\cite{StrichartzTransformations}.  The appropriate substitution operation preserves spectral similarity of graph Laplacians according to the following result from~\cite{MT}

\begin{lemma}[\protect{Lemma~4.7 in~\cite{MT}}]\label{lem:glue}
If\/ $\Gamma$ is a graph with vertex set $V$ and is symmetric with respect to $\tilde{V}\subset V$ in the sense that every bijection of $\tilde V$ extends to a graph automorphism of $\Gamma$ then the graph Laplacian $\Delta^\Gamma$ of\/ $\Gamma$ is spectrally similar to $\Delta^{\tilde{\Gamma}}$ for the complete graph $\tilde \Gamma$ on $\tilde V$, with some functions $\phi$ and $\tilde\phi$.  Now take disjoint copies of $\Gamma$ over some index set $\mathcal{A}$ and form a quotient graph $\Gamma'$ by identifying vertices in $V\times\mathcal{A}$ via an equivalence relation on $\tilde V\times \mathcal{A}$ for which equivalence classes are finite. Repeat the same quotient construction on $\tilde\Gamma\times\mathcal{A}$ to get a quotient graph $\tilde \Gamma'$.  Then $\Delta^{\Gamma'}$ is spectrally similar to $\Delta^{\tilde \Gamma'}$ with the same functions $\phi$ and $\tilde\phi$.
\end{lemma}

\begin{proposition}\label{prop:Gndec}
The Laplacian $\Delta^J_n$ on $J_n$ is spectrally similar to $\Delta^G_n$ on $G_n$ with exceptional set $E=\{-1\}$ and decimation function $R_3(z)=\frac32(z^2+2z)$. If $z\neq-1$ and $R_3(z)$ is an eigenvalue of $\Delta^G_n$ with eigenfunction $\tilde f$ then the corresponding eigenfunction of $\Delta^J_n$ with eigenvalue $z$ is obtained by extending $\tilde f$ to each point $x\in V^J_n\setminus V^G_n = V^H_n$ so that $f(x)=\frac1{3(1+z)} \sum_{y\sim_{J_n} x} \tilde f(y)$.  Moreover, $-1$ is an eigenvalue of $\Delta^J_n$ with a $\frac12(3^n+3)$ dimensional eigenspace consisting of functions $g$ that vanish at each $x\in V^H_n\subset V^J_n$ and such that $\sum_{y\sim_{J_n} x} g(y)=0$ for all such $x$.
\end{proposition}
\begin{proof}
Recall we wrote $F_w=F_{w_1}\circ\dotsc F_{w_n}$ where $w_1\dotsm w_n$ is a length $n$ word with letters $w_j\in\{1,2,3\}$. 
The vertices $V^J_n$ are then all points of the form $F_w(p_j)$  as $w$ ranges over $n$-words and $j\in\{0,1,2,3\}$, and the inductive  definition of the edges of $J_n$ implies there is an edge between $F_w(p_j)$ and $F_{w'}(p_k)$ if and only if $w=w'$ and $j\neq k$.  This exhibits $J_n$ as a quotient as in Lemma~\ref{lem:glue} of copies of $J_0$ indexed by the length $n$-words $\{1,2,3\}^n$.

Adopting the notation of the lemma, $\Gamma=J_0$, $\mathcal{A}=\{1,2,3\}^n$, and the equivalence relation is that $(p_j,w)$ and $(p_k,w')$ are equivalent when $F_w(p_j)=F_{w'}(p_k)$.  Observe that this latter can never occur if $j$ or $k$ is zero, so the equivalence relation is really on $\tilde{V}\times\mathcal{A}$, where $\tilde{V}=\{p_1,p_2,p_3\}=V^G_0$, and it is easily seen that $V^J_0$ is symmetric with respect to $V^G_0$. 
According to  Lemma~\ref{lem:glue}, then $\Gamma$ is spectrally similar to the complete graph $\tilde\Gamma$ on $\tilde V$, which we observe is actually $G_0$.  In fact this is the case $N=3$ of the spectral similarity in Example~\ref{ex:starmap}, so we know the functions $\phi_3$ and $\tilde\phi_3$ explicitly.

We have seen that $J_n=\Gamma'$, the quotient in  Lemma~\ref{lem:glue}. At the same time, the quotient $\tilde\Gamma'$ is that obtained by gluing the corresponding copies of $\tilde\Gamma=V^G_0$, and the same reasoning as used on $J_n$ shows the result is $\tilde\Gamma'=G_n$.  Now the result of the lemma says that $J_n$ and $G_n$ are spectrally similar with functions $\phi(z)=\phi_3(z)= 1+z-\frac1{1+z}$ and $\tilde\phi(z)=\tilde\phi_3(z) = \frac2{3(1+z)}$, giving the formula for $R_3(z)$.   Moreover, the discussion in Example~\ref{ex:starmap} tells us that if $z\neq-1$ and $R_3(z)\in\sigma(\Delta^G_n)$ with eigenfunction $\tilde f$ then we can extend to an eigenfunction $f$ of $\Delta^J_n$ with eigenvalue $z$ by taking $f(F_w(p_0))=\frac1{3(1+z)}\sum_{j=1}^3 \tilde f(F_w(p_j))$.  At the same time, $z=-1$ is an eigenvalue of $\Delta^J_n$  with eigenspace consisting  of functions that vanish at all points $F_w(p_0)$ and have mean zero over the neighbors of such points. Evidently this space involves the values on the $\frac12(3^{n+1}+3)$ points in $V^G_n$ subject to $3^n$ constraints (because the mean is zero on each of the $3^n$ cells), so is dimension $\frac12(3^n+3)$.\end{proof}

\begin{proposition}\label{prop:Hndec}
The Laplacian $\Delta^J_n$ on $J_n$ is spectrally similar to $\Delta^H_n$ on $H_n$ with exceptional set $E=\{-1\}$ and decimation function $R_2(z)=2(z^2+2z)$. If $z\neq-1$ and $R_2(z)$ is an eigenvalue of $\Delta^H_n$ with eigenfunction $\tilde f$ then the corresponding eigenfunction of $\Delta^J_n$ with eigenvalue $z$ is obtained by extending $\tilde f$ to each point $x\in V^J_n\setminus V^G_n = V^H_n$ so that $f(x)=\frac1{2(1+z)} \sum_{y\sim_{J_n} x} \tilde f(y)$.
\end{proposition}
\begin{proof}
We again use Lemma~\ref{lem:glue}, in particular following its notation.  This time $\Gamma$ is the graph called $\Gamma_2$ in Example~\ref{ex:starmap}, so the vertex set $V$ has three elements which we label $\{q_0,q_1,q_2\}$, with $q_0$ being the one incident on two edges. Set $\tilde V=\{q_1,q_2\}$. Evidently $V$ is symmetric with respect to $\tilde V$ and the complete graph $\tilde\Gamma$ on $\tilde V$ is a single edge.  The spectral symmetry functions are $\phi_2$ and $\tilde\phi_2$ from Example~\ref{ex:starmap}.

Our index set $\mathbb{A}$ will be the edge set of $H_n$, including the loops at the boundary vertices.  The disjoint union of copies of $\tilde\Gamma$ indexed by $\mathcal{A}$ has vertex set $\tilde V\times\mathcal{A}$. On it we define an equivalence relation such that the quotient $\tilde\Gamma'$ is precisely $H_n$. This is easy to do: whenever $(q_1,e), (q_2,e)$ are the two points in $\tilde V\times\{e\}$ for some $e\in\mathcal{A}$ we simply identify them with the vertices of $V^H_n$ that are the endpoints of $e$. The identification is a bijection unless $e$ is a loop, in which case it is two-to-one.  Evidently this is a surjection $T:\tilde V\times\mathcal{A}$ to $V^H_n$, so we obtain an equivalence relation $(q_1,e_1)\sim(q_2,e_2)$ if  $T(q_1,e_1)=T(q_2,e_2)$.  This equivalence is not uniquely defined because we must choose one of the two possible bijections at each $e$, except in the case that $e$ is a loop. However, every edge in the quotient comes from an edge in a copy of $\tilde\Gamma$, so for any set of choices of the bijections we see that there is an edge between the equivalence classes of $(q_1,e_1)$ and $(q_2,e_2)$ in $\tilde\Gamma'$ if and only if we can choose representatives such that $e_1=e_2=e$, which happens precisely when $T(q_1,e)$ and $T(q_2,e)$ are joined by the edge $e$ in $H_n$. Moreoever this is a loop exactly when $(q_1,e)\sim(q_2,e)$, meaning $T(q_1,e)=T(q_2,e)$, which occurs exactly when $e$ is a loop in $H_n$. Thus $\tilde\Gamma'=H_n$.

Now we wish to identify the quotient $\Gamma'$.  We had $T:\tilde V\times\mathcal{A}\to V^H_n$ and we extend it to $V\times\mathcal{A}\to V^J_n$ by defining it on $(V\setminus\tilde V)\times\mathcal{A}=\{q_0\}\times\mathcal{A}$ as follows. Given $(q_0,e)$ we have defined $T(q_1,e)$ and $T(q_2,e)$ to be endpoints of the edge $e$ from $H_n$. If $e$ is a loop then $T(q_1,e)=T(q_2,e)\in V^H_n\subset V^J_n$ and this vertex is joined in $J_n$ to a unique vertex $p\in V^G_0\subset V^J_n$; in this case let $T(q_0,e)=p$.  Evidently this extension is a bijection from $\{(q_0,e): e\text{ is a loop in $H_n$}\}$ to $V^G_0$.

If $e$ is an edge in $H_n$ with endpoints  $T(q_1,e)=F_w(p_0)$ and $T(q_2,e)=F_{w'}(p_0)$ from $V^H_n$, where $w,w'$ are words of length $n$, then our definition of $H_n$ says there is a unique point $x\in V^G_n$ so that $F_w(p_0)$ and $F_{w'}(p_0)$ are both joined to $x$ by edges in $J_n$. We define $T(q_0,e)=x$. Since every point of $V^G_n\setminus V^G_0$ is of the form $F_w(p_j)=F_{w'}(p_k)$ for some length $n$ words $w\neq w'$ and $j\neq k\in\{1,2,3\}$ we see there is a corresponding edge of $H_n$ from $F_w(p_0)$ to $F_{w'}(p_0)$ and thus $T$ is bijective from $\{(q_0,e): e\text{ is an edge of $H_n$}\}\subset V \times\mathcal{A}$ to $V^G_n\setminus V^G_0$.

Our extension is a surjection $T:V\times\mathcal{A}\to V^J_n$ that is bijective from $(V\setminus \tilde V)\times\mathcal{A}$ to $V^G_n$. Taking the quotient $\Gamma'$ by the previously defined equivalence on $\tilde V\times\mathcal{A}$ we see that $T$ identifies the vertex set of $\Gamma'$ with $V^J_n$. To complete the proof we need to show that pairs of equivalence classes from $\Gamma'$ that are joined by an edge correspond under $T$ to pairs of points in $V^J_n$ that are joined by an edge of $J_n$. To this end, observe that edges in $V\times\mathcal{A}$ are between pairs $(q_0,e)$ and $(q_j,e)$ for $j=1,2$ so edges in $\Gamma'$ are between pairs of equivalence classes represented by $(q_0,e)$ and $(q_j,e)$, $j=1,2$.  We will deal with the case where $e$ is a loop of $H_n$ separately from the case where it is an edge.

Fix equivalence classes in $\Gamma'$ that are joined by an edge and choose representatives as described above.  Bijectivity of $T:(V\setminus \tilde V)\times\mathcal{A}\to V^G_n$ says the equivalence class of $(q_0,e)$ is a one-element set corresponding to the point $x=T(q_0,e)\in V^G_n$. If $e$ is an edge (not a loop) of $H_n$ then $x\notin V^G_0$ so $x=F_w(p_j)=F_{w'}(p_k)$ for exactly two distinct length $n$ words $w\neq w'$ and corresponding $j\neq k\in\{1,2,3\}$. Our construction ensures $e$ is the edge of $H_n$ between $F_w(p_0)$ and $F_{w'}(p_0)$ and thus  $T(q_j,e)$ is one of $F_w(p_0)$ or $F_{w'}(p_0)$.  However, the edges  in $J_n$ that are incident at $x$ connect to precisely the points $F_w(p_0)$ and $F_{w'}(p_0)$, so we conclude that $T(q_0,e)$ and $T(q_j,e)$ are connected by an edge in $J_n$. A similar argument applies if $e$ is a loop of $H_n$, as then $x\in V^G_0$ is $p_k=F_{k^n}(p_k)$ and $T(q_j,e)=F_{k^n}(p_0)$ for both $j=1,2$. This latter is exactly the neighbor of $p_k$ in $J_n$.

Conversely, fix an edge of $J_n$ which connects $F_w(p_0)$ to $x\in V^G_n$. If $x\notin V_0$ we can write $x=F_w(p_k)=F_{w'}(p_l)$ for some length $n$ words $w\neq w'$ and some $k\neq l\in\{1,2,3\}$. The definition of $H_n$ gives an edge $e$ between $F_w(p_0)$ and $F_{w'}(p_0)$. Our construction provides $T(q_0,e)=x$ and $\{T(q_j,e):j=1,2\}=\{F_w(p_0),F_{w'}(p_0)\}$, so $T^{-1}(x)=\{(q_0,e)\}$ and  $T^{-1}(F_w(p_0))$ is the equivalence class of either $(q_1,e)$ or $(q_2,e)$, both of which neighbor $\{(q_0,e)\}$ in $\Gamma'$. If, on the other hand, $x\in V^G_0$ then $x=p_k$ neighbors $F_{k^n}(p_0)$ in $J_n$, and the definition of $T$ gives $T^{-1}(x)=(q_0,e)$ for the loop of $H_n$ attached at $F_{k^n}(p_0)$. Thus $(q_j,e)\in T^{-1}(F_{k^n}(p_0))$ for both $j=1,2$ and this equivalence class is a neighbor of $T^{-1}(x)$ in $\Gamma'$.

The preceding proves that $T$ induces a graph isomorphism between $\Gamma'$ and $J_n$, so we can apply Lemma~\ref{lem:glue} to conclude that $J_n$ is spectrally similar to $H_n$ with functions (from Example~\ref{ex:starmap}) $\phi=\phi_2=1+z-\frac1{1+z}$  and $\tilde\phi_2=\frac1{2(1+z)}$.  The stated expression for $R_2(z)$ follows, as does the fact that the exceptional set is $E=\{-1\}$. Moreover, we find from the discussion in Example~\ref{ex:starmap} that if $z\neq-1$ and $\tilde f$ is an eigenfunction of $\Delta^H_n$ with eigenvalue $R_2(z)$ then we can extend $\tilde f$ to an eigenfunction $f$ of $\Delta^J_n$ with eigenvalue $z$ by setting $f(x)=\frac1{2(1+z)}\sum_{y\sim_{J_n}x}\tilde f(y)$  at points of $V^G_n\setminus V^G_0$, while if $x\in V^G_0$ and $y\in V^H_n\subset V^J_n$ is its unique neighbor in $J_n$ then $f(x)=\tilde f(y)$.

The discussion in Example~\ref{ex:starmap} also suggests that we could construct  eigenfunctions $g$ of $\Delta^J_n$ with eigenvalue $-1$ by requiring that  $g$ vanish at each $x\in V^G_n\subset V^J_n$ and have $\sum_{y\sim_{J_n}x}g(y)=0$ for all such $x$. However these requirements imply $g$ is identically zero, so there are no such eigenfunctions. To see this, take points $y,y;\in V^H_n$ and let $x\in V^G_n$ be their common neighbor in $J_n$.  The constraint equation says that if $g(y)=0$ then also $g(y')=0$, so by the obvious path connectedness of $J_n$ we determine that if $g$ satisfying the constraint equation vanishes at some $y\in V^H_n$ then it vanishes on all of $V^H_n$. Since the other constraint says $g$ vanishes on $V^G_n$, we conclude that if $g$ vanishes at some $y\in V^H_n$ it is identically zero on $J_n$.  However, if we take $x=p_1$ in the constraint equation then the fact that $p_1$ has only one neighbor in $J_n$ ensures $g$ will vanish at that neighbor, which is a point of $V^H_n$.
\end{proof}

\begin{remark}
There is a subtlety about the application of Lemma~\ref{lem:glue} in the proof of Theorem~\ref{prop:Hndec} that may warrant additional explanation. It is clear that the identification of two endpoints of an edge to a single point gives a loop, but when we identify the two endpoints $q_1,q_2$ of a copy of $\Gamma=\Gamma_2$ to the same single point it appears to give a double edge rather than the single edge we expect to see in $J_n$. The resolution of this apparent difficulty is embedded in the proof of Lemma~\ref{lem:glue}, so is not visible in the proof of Theorem~\ref{prop:Hndec}. However, we can see what happens by looking at the weights for the connections from $p_j\in V^G_0$ to the adjacent point $y=F_{j^n}(p_0)$. In $J_n$ we have a single edge to $y$ with weight $1=1/\deg(p_j)$. In $\Gamma'$ our construction appears to provide two edges from $p_j$ to $y$, but since each has adjacency weight $1/\deg(p_j)=1/2$ this is equivalent to having a single edge of weight $1$.
\end{remark}

\begin{theorem}\label{thm:main}
The sets $ \sigma(\Delta^G_n)\setminus\{-\frac32\}$ and $\sigma(\Delta^H_n)\setminus\{-2\}$ are bijective under the map $\zeta\mapsto\frac43\zeta$. Moreover, if we define for $f:V^G_n\to\mathbb{C}$ and $g:V^H_n\to\mathbb{C}$
\begin{align*}
	\Phi_3f(y) &= \frac13 \sum_{x\sim_{J_n}y}f(x) \quad \text{ for $x\in V^H_n$, and}\\
	\Phi_2g(x) &=
		\begin{cases}   \frac12\sum_{y\sim_{J_n}x}g( y) & \text{ if $x\in V^G_n\setminus V^G_0$}\\
					\sum_{y\sim_{J_n}x}g( y)  & \text{ if $x\in V^G_0$}
					\end{cases}
	\end{align*}
then $\Phi_3$ is a bijection from the eigenspace of $\Delta^G_n$ with eigenvalue $\zeta\in \sigma(\Delta^G_n)\setminus\{-\frac32\}$ to the eigenspace of $\Delta^H_n$ with eigenvalue $\frac43\zeta\in \sigma(\Delta^H_n)\setminus\{-2\}$, and $\Phi_2$ is its inverse.
\end{theorem}
\begin{proof}
Propositions~\ref{prop:Gndec} and~\ref{prop:Hndec} tell us that the following maps are surjections
\begin{align*}
	R_3: \sigma(\Delta^J_n)\setminus E\to\sigma(\Delta^G_n)\setminus\{-\frac32\}\\
	R_2: \sigma(\Delta^J_n)\setminus E\to\sigma(\Delta^H_n)\setminus\{-2\}
	\end{align*}
and it is easily seen that $-1$ is the only critical point of both $R_3$ and $R_2$, so these maps are actually two-to-one covers. By inspection of the formulas for $R_3$ and $R_2$ we conclude that
\begin{equation*}
	 \sigma(\Delta^G_n)\setminus\{-\frac32\} \xrightarrow{\zeta\mapsto\frac43\zeta} \sigma(\Delta^H_n)\setminus\{-2\}
	 \end{equation*}
is a bijection. 

Now comparing the definition of  the maps $\Phi_2$ and $\Phi_3$ with the extensions given in Propositions~\ref{prop:Gndec} and~\ref{prop:Hndec} we see that if $z\in\sigma(\Delta^J_n)\setminus E$ with eigenfunction $f$ then $f|_{V^G_n}=(1+z)^{-1} \Phi_2 \bigl( f|_{V^H_n})$ and $f|_{V^H_n}=(1+z)^{-1} \Phi_3 \bigl( f|_{V^G_n})$. Since $f|_{V^G_n}$ is an eigenfunction of $\Delta^G_n$ with eigenvalue $R_3(z)$ and $f|_{V^H_n}$ is an eigenfunction of $\Delta^H_n$ with eigenvalue $R_2(z)=\frac43R_3(z)$, we conclude that $\Phi_3$ maps the eigenspace of $\Delta^G_n$ with eigenvalue $\zeta\in\sigma(\Delta^G_n)\setminus\{-\frac32\}$ bijectively to the eigenspace of $\Delta^H_n$ with eigenvalue $\frac43$. Repeating the argument for $\Phi_2$ shows that it is inverse to $\Phi_3$ on these spaces. 
This last point may seem counterintuitive; it is more readily apparent when one recalls that these eigenfunctions on $J_n$ are orthogonal to the functions having mean zero on all sets $\{F_w(p_j):j=1,2,3\}$ where $w$ is a length $n$ word, because the latter are eigenfunctions with eigenvalue $-1$.
\end{proof}

\section{Spectral self-similarity of $\Delta^G_n$ and $\Delta^H_n$}\label{sec:specselfsim}
Theorem~\ref{thm:main} allows us to get the spectrum of $H_n$ from that of $G_n$ or vice-versa.  The spectrum for the graphs $G_n$ has been known for a long time~\cite{RammalToulouse, FukushimaShima} as a consequence of the spectral self-similarity of the Laplacians $\Delta^G_n$ given in Theorem~\ref{thm:SGspect} below.   While the spectrum for $H_n$ was given in~\cite{GrigorchukSunic}, we obtain it by transferring this spectral self-similarity to the Laplacians $\Delta^H_n$ via Theorem~\ref{thm:main}.

\begin{theorem}[\protect{\cite{RammalToulouse, FukushimaShima,DSV}}]\label{thm:SGspect}
For $n\geq0$ the Laplacian $\Delta^G_{n+1}$ is spectrally similar to $\Delta^G_n$ with exceptional set $E_G=\{-\frac32, -\frac54,-\frac12\}$ and decimation function $R_G(z)=z(4z+5)$. If $z\notin E_G$ and $R_G(z)$ is an eigenvalue of $\Delta^G_n$ with eigenfunction $\tilde f$ then the corresponding eigenfunction of $\Delta^G_{n+1}$ with eigenvalue $z$ is obtained by extending $\tilde f$ to each point of $V^G_{n+1}\setminus V^G_n$ in the following manner: for $w$ a length $n$ word we have $\tilde f$ at  $F_w(p_j)$, $j=1,2,3$;  write $m_j=F_{j+1}(p_{j-1})$ for the vertex opposite $p_j$ in $V_1$ (all subscripts are taken modulo $3$ so as to lie in $\{1,2,3\}$), and set
\begin{equation}\label{eq:Gefnextension}
	f(F_w(m_j))=\frac{ f(F_w(p_j)) + 2(z+1)\bigl( f(F_w(p_{j+1}))+f(F_w(p_{j-1}))\bigr)}{(4z+5)(2z+1)}.
	\end{equation}
\end{theorem}

\begin{remark}
The interested reader may care to prove this theorem herself by checking that $D=\Delta^G_1$ is spectrally similar to $\tilde D=\Delta^G_0$ with $\phi=z(2z+3)/(2z+1)$ and $\tilde\phi=(2z+3)/(4z+5)(2z+1)$ and exceptional set $E_G$, then working as in Proposition~\ref{prop:Gndec} to write $G_n$ as a quotient of copies of $G_0$ indexed by words in $\{1,2,3\}^n$ and checking that the same quotient of copies of $G_1$ is $G_{n+1}$. Alternatively, an elementary and very readable treatment appears in~\cite[Chapter 3]{Strichartzbook}, with a different normalization for the Laplacian.
\end{remark}

The bijection in Theorem~\ref{thm:main} allows us to transfer this spectral similarity to the graphs $H_n$.
\begin{theorem}\label{thm:Hndecimation}
The Laplacian $\Delta^H_{n+1}$ is spectrally similar to $\Delta^H_n$ with exceptional set $E_H=\{-2,-\frac53,-1,-\frac23\}$ and decimation function $R_H(z)= z(3z+5)$. If $z\notin E_H$ and $R_H(z)$ is an eigenvalue of $\Delta^H_n$ with eigenfunction $f$ then the corresponding eigenfunction $g$ of $\Delta^H_{n+1}$ with eigenvalue $z$ is obtained as follows. A point $x\in V^H_{n+1}$ is $F_{wj}(p_0)$ for some $w$ of length $n$ and $j\in\{1,2,3\}$. We have $y_0=F_w(p_0)\in V^H_n$; write $y_k$, $k=1,2,3$ for the neighbors of $y_0$ in $H_n$, choosing the labels so that $y_0$ is joined to $y_k$ via $F_w(p_k)$ in $J_n$. (Note that $y_k=y_0$ if $F_w(p_k)\in V^G_0$, meaning there is a loop at $y_0$.) Then
\begin{equation}\label{eq:Hnextension}
	g(x) =   \frac{2z+3}{6(4z+5)(2z+1)}\Bigl( (4z+5)f(y_0) +(4z+3)f(y_j) + f(y_{j+1})+f(y(j-1))\Bigr) 
	\end{equation}
where, as usual, $j\pm1$ are taken modulo $3$ so as to give values in $\{1,2,3\}$.
\end{theorem}

\begin{figure}
\begin{tikzpicture}[scale=1.3]
\draw  (0,0) node (11) {} -- (1.25,0) node (12) {} -- (0.625,1.0825317547305) node (13) {}--cycle;
\fill  (0,0) circle[radius=1.3pt]  (1.25,0)  circle[radius=1.3pt]  (0.625,1.0825317547305) circle[radius=1.3pt] ;
\draw  (1.25,0) node (21) {} -- (2.5,0) node (22){} -- (1.875,1.0825317547305) node (23) {}--cycle;
\fill  (1.25,0) circle[radius=1.3pt]  (2.5,0)  circle[radius=1.3pt]  (1.875,1.0825317547305) circle[radius=1.3pt] ;
\draw  (0.625,1.0825317547305) node (31) {} -- (1.875,1.0825317547305) node (32) {} -- (1.25,2.1650635094611) node (33) {}--cycle;
\fill  (0.625,1.0825317547305) circle[radius=1.3pt]  (1.875,1.0825317547305)  circle[radius=1.3pt]  (1.25,2.1650635094611) circle[radius=1.3pt] ;
\draw  (2.5,0) node (41) {} -- (3.75,0) node (42) {} -- (3.125,1.0825317547305) node (43) {}--cycle;
\fill  (2.5,0) circle[radius=1.3pt]  (3.75,0)  circle[radius=1.3pt]  (3.125,1.0825317547305) circle[radius=1.3pt] ;
\path  (0.625,0.36084391824352) node (10) {b} (1.875,0.36084391824352) node (20) {a}  (1.25,1.4433756729741) node (30) {c}  (3.125,0.36084391824352) node (40) {d};
\draw[dashed] (10)--(20)--(30)--(10) (20)--(40); 
\path (12) node[below] {$\scriptstyle\frac12(a+b)$} (23) node[right] {$\scriptstyle\frac12(a+c)$} (22) node[below] {$\scriptstyle\frac12(a+d)$};
\draw[dotted] (10)--(11) (10)--(13) (10)--(12) (20)--(21) (20)--(22) (20)--(23) (30)--(31) (30)--(32) (30)--(33) (40)--(41) (40)--(42) (40)--(43); 
\end{tikzpicture}\ \ \ \
\begin{tikzpicture}[scale=1.3]
\draw  (0,0) node (11) {} -- (1.25,0) node (12) {} -- (0.625,1.0825317547305) node (13) {}--cycle;
\fill  (0,0) circle[radius=1.3pt]  (1.25,0)  circle[radius=1.3pt]  (0.625,1.0825317547305) circle[radius=1.3pt] ;
\draw  (1.25,0) node (21) {} -- (2.5,0) node (22) {} -- (1.875,1.0825317547305) node (23) {}--cycle;
\fill  (1.25,0) circle[radius=1.3pt]  (2.5,0)  circle[radius=1.3pt]  (1.875,1.0825317547305) circle[radius=1.3pt] ;
\draw  (0.625,1.0825317547305) node (31) {} -- (1.875,1.0825317547305) node (32) {} -- (1.25,2.1650635094611) node (33) {}--cycle;
\fill  (0.625,1.0825317547305) circle[radius=1.3pt]  (1.875,1.0825317547305)  circle[radius=1.3pt]  (1.25,2.1650635094611) circle[radius=1.3pt] ;
\draw  (2.5,0) node (41) {} -- (3.75,0) node (42) {} -- (3.125,1.0825317547305) node (43) {}--cycle; \fill  (2.5,0) circle[radius=1.3pt]  (3.75,0)  circle[radius=1.3pt]  (3.125,1.0825317547305) circle[radius=1.3pt] ;
\draw  (0,2.1650635094611) node (51) {} -- (1.25,2.1650635094611) node (52) {} -- (0.625,3.2475952641916) node (53) {}--cycle; \fill   (0,2.1650635094611) circle[radius=1.3pt] (1.25,2.1650635094611)  circle[radius=1.3pt]  (0.625,3.2475952641916) circle[radius=1.3pt] ;
\draw (-1.25,0) node (61) {} -- (0,0) node (62) {} -- (-.625,1.0825317547305) node (63) {}--cycle;
\fill   (-1.25,0)  circle[radius=1.3pt]  (-0.625,1.0825317547305) circle[radius=1.3pt] ;
\path (11) node[below] {$\scriptstyle\frac12(a+b)$} (12) node[below] {$\alpha$} (13) node[left]{$\beta$} (22) node[below] {$\scriptstyle\frac12(a+d)$} (23) node[right]{$\gamma$} (33) node[right]{$\scriptstyle\frac12(a+c)$};
\path  (0.625,0.36084391824352) node (10) {s} (1.875,0.36084391824352) node (20) {t}  (1.25,1.4433756729741) node (30) {u} (3.125,0.360843918243) node[coordinate] (40) {} (0.625,2.52590742770462) node[coordinate] (50) {}    (-0.625,0.36084391824352) node[coordinate] (60) {} ;
\draw[dotted] (10)--(11) (10)--(13) (10)--(12) (20)--(21) (20)--(22) (20)--(23) (30)--(31) (30)--(32) (30)--(33) (40)--(41) (40)--(42) (40)--(43) (50)--(51) (50)--(53) (50)--(52) (60)--(61) (60)--(63) (60)--(62); 
\draw[dashed] (60)--(10)--(20)--(30)--(10) (20)--(40) (30)--(50) ; 
\end{tikzpicture}
\caption{\protect{Computation of~\eqref{eq:Hnextension}: values $a,b,c,d$ on $H_n$ (dashed lines) are extended to $G_n$ (solid lines) via $\Phi_2$, then to $\alpha,\beta,\gamma$ in $G_{n+1}$ by~\eqref{eq:Gefnextension} and thus to $s,t,u$ on $H_{n+1}$ via $\Phi_3$. Dotted lines show $J_n$ and $J_{n+1}$.}}\label{fig:computeHnefnextend}
\end{figure}
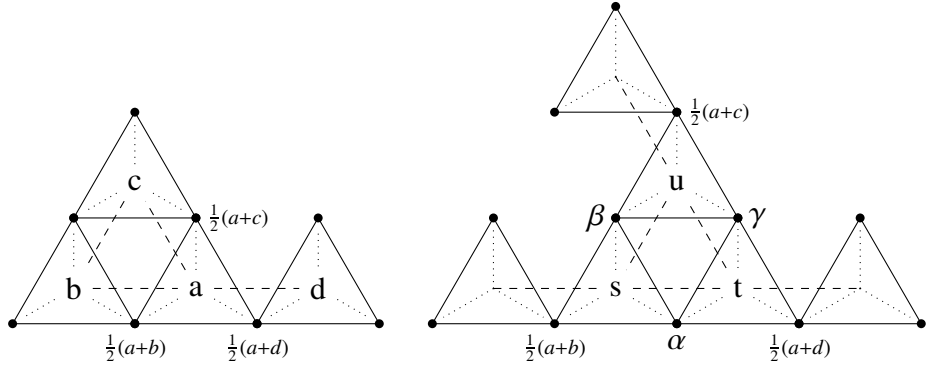

\begin{proof}
From Theorem~\ref{thm:main} we have a bijection between $\sigma(\Delta^G_m)\setminus\{-\frac32\}$ and $\sigma(\Delta^H_m)\setminus\{-2\}$, simply by $\zeta\mapsto\frac43\zeta$. Suppose $z\in\sigma(\Delta^H_{n+1})\setminus E_H$. Then this bijection gives $\frac34z\in\sigma(\Delta^G_m)\setminus (E_G\cup\{-\frac34\})$, and therefore Theorem~\ref{thm:SGspect} says that $R_G(\frac34z)\in\sigma(\Delta^G_n)$.  Moreover we have $R_G^{-1}(-\frac32)=\{-\frac34,-\frac12\}\subset(E_G\cup\{-\frac34\})$, so in fact $R_G(\frac34z)\in\sigma(\Delta^G_n)\setminus\{-\frac32\}$ and we can again use the bijection from Theorem~\ref{thm:main}, this time with $m=n$, to get that $\frac43R_G(\frac34z)\in\sigma(\Delta^H_n)\setminus\{-2\}$.  Hence $R_H(z)=\frac43R_G(\frac34z)= z(3z+5)$ is a spectral decimation map from $\Delta^H_{n+1}$ to $\Delta^H_n$ with the stated exceptional set.

Given an eigenfunction as in the statement of the theorem we apply $\Phi_2$ to get an eigenfunction of $\Delta^G_n$, then the map in~\eqref{eq:Gefnextension}, and finally $\Phi_3$ to get an eigenfunction of $\Delta^H_{n+1}$. For the reader's convenience we illustrate with the sequence of diagrams in Figure~\ref{fig:computeHnefnextend}. One is given the values $a,b,c,d$ on a piece of $H_n$ shown on the left (dashed lines). Treating these as values on $J_n$ (dotted lines), apply $\Phi_2$ to obtain values at the vertices in $G_n$. Then~\eqref{eq:Gefnextension} gives
\begin{align*}
	2\alpha	&= \mu(2a+b+d)) + \nu(a+c)\\
	2\beta  	&= \mu(2a+b+c)) + \nu(a+d)\\
	2\gamma	&= \mu(2a+c+d)) + \nu(a+b)
	\end{align*}
where $\mu=\frac{2(z+1)}{(4z+5)(2z+1)}$ and $\nu=\frac1{(4z+5)(2z+1)}$. Apply $\Phi_3$ to get
\begin{align}
	6t&=2\alpha+2\gamma+a+d =  \mu(4a+b+c+2d)) + \nu(2a+b+c) +a+d \notag\\
	&= \frac{2z+3}{(4z+5)(2z+1)}\Bigl( (4z+5)a + (b+c)+ (4z+3)d\Bigr) \label{eq:Hndecintermedcomputation1}
	\end{align}
This is precisely~\eqref{eq:Hnextension}. To see this, write $t=g(F_{wj}(p_0))$. Then $f(F_w(p_0))=a$ is apparent, and $y_0=F_w(p_0)$. There are edges in $H_n$ from $y_0$ to each $y_k$, and these correspond to paths in $J_n$ via points $F_w(p_k)\in V^G_n$. Since $d$ has a different weight to $b$ and $c$ in~\eqref{eq:Hndecintermedcomputation1} it is important that the edge between the corresponding vertices of $H_n$ goes via $F_w(p_j)$, and thus $d=f(y_j)$.  The other two values $f(y_k)$, $k\in\{1,2,3\}\setminus\{j\}$ are $b$ and $c$.
\end{proof}

\section{Eigenfunctions of $G_n$ and $H_n$}\label{sec:efns}

It is well-known that one can give a complete description of all eigenfunctions of $\Delta^G_n$ using the results in the previous section. This is because eigenfunctions on $G_{n+1}$ either have eigenvalue in $E_G$ or are obtained by decimation from Theorem~\ref{thm:SGspect}: one simply makes a direct construction of eigenfunctions  with eigenvalues in $E_G$ and counts to ensure they generate the whole spectrum. In this section we use the same idea to construct the eigenfunctions of $\Delta^H_n$ and determine multiplicities in the spectrum.  A convenient approach is to state known results  for $\Delta^G_n$ as propositions and derive the implications for $\Delta^H_n$ sequentially.

\begin{proposition}[\protect{\cite{DSV}, see also~\cite[Chapter 3]{Strichartzbook}}]\label{prop:Gnexcept1}
The exceptional value $-\frac32$ is an eigenvalue of $G_n$ for all $n\geq0$; the corresponding eigenfunctions are exactly those functions satisfying $\sum_{j=1}^3 f(F_w(p_j))=0$ for all words $w$ of length $n$ and form a $\frac12(3^n+3)$ dimensional eigenspace. 
\end{proposition}

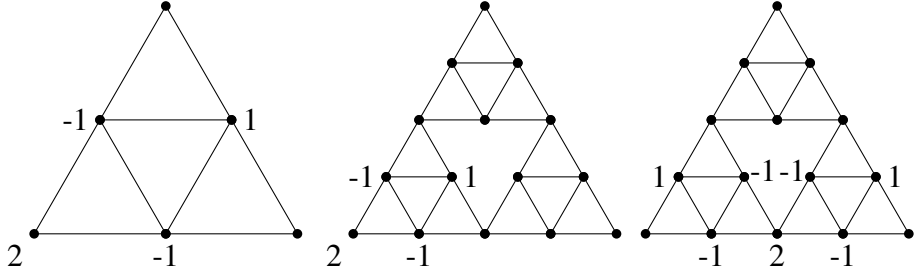
\begin{figure}
\begin{tikzpicture}[scale=1.4]
\draw  (0,0) node (11) {} -- (1.25,0) node (12) {} -- (0.625,1.0825317547305) node (13) {}--cycle;\fill  (0,0) circle[radius=1.3pt]  (1.25,0)  circle[radius=1.3pt]  (0.625,1.0825317547305) circle[radius=1.3pt] ;\draw  (1.25,0) node (21) {} -- (2.5,0) node (22) {} -- (1.875,1.0825317547305) node (23) {}--cycle;\fill  (1.25,0) circle[radius=1.3pt]  (2.5,0)  circle[radius=1.3pt]  (1.875,1.0825317547305) circle[radius=1.3pt] ;\draw  (0.625,1.0825317547305) node (31) {} -- (1.875,1.0825317547305) node (32) {} -- (1.25,2.1650635094611) node (33) {}--cycle;\fill  (0.625,1.0825317547305) circle[radius=1.3pt]  (1.875,1.0825317547305)  circle[radius=1.3pt]  (1.25,2.1650635094611) circle[radius=1.3pt] ;
\path (11) node[below left] {2} (12)  node[below]  {-1} (13) node[left] {-1}  (23) node[right] {1};
\end{tikzpicture}
\begin{tikzpicture}[scale=1.4]
\draw  (0,0) node (111) {} -- (0.625,0) node (112) {} -- (0.3125,0.54126587736527) node (113) {}--cycle;\fill  (0,0) circle[radius=1.3pt]  (0.625,0)  circle[radius=1.3pt]  (0.3125,0.54126587736527) circle[radius=1.3pt] ;\draw  (0.625,0) node (121) {} -- (1.25,0) node (122) {} -- (0.9375,0.54126587736527) node (123) {}--cycle;\fill  (0.625,0) circle[radius=1.3pt]  (1.25,0)  circle[radius=1.3pt]  (0.9375,0.54126587736527) circle[radius=1.3pt] ;\draw  (0.3125,0.54126587736527) node (131) {} -- (0.9375,0.54126587736527) node (132) {} -- (0.625,1.0825317547305) node (133) {}--cycle;\fill  (0.3125,0.54126587736527) circle[radius=1.3pt]  (0.9375,0.54126587736527)  circle[radius=1.3pt]  (0.625,1.0825317547305) circle[radius=1.3pt] ;\draw  (1.25,0) node (211) {} -- (1.875,0) node (212) {} -- (1.5625,0.54126587736527) node (213) {}--cycle;\fill  (1.25,0) circle[radius=1.3pt]  (1.875,0)  circle[radius=1.3pt]  (1.5625,0.54126587736527) circle[radius=1.3pt] ;\draw  (1.875,0) node (221) {} -- (2.5,0) node (222) {} -- (2.1875,0.54126587736527) node (223) {}--cycle;\fill  (1.875,0) circle[radius=1.3pt]  (2.5,0)  circle[radius=1.3pt]  (2.1875,0.54126587736527) circle[radius=1.3pt] ;\draw  (1.5625,0.54126587736527) node (231) {} -- (2.1875,0.54126587736527) node (232) {} -- (1.875,1.0825317547305) node (233) {}--cycle;\fill  (1.5625,0.54126587736527) circle[radius=1.3pt]  (2.1875,0.54126587736527)  circle[radius=1.3pt]  (1.875,1.0825317547305) circle[radius=1.3pt] ;\draw  (0.625,1.0825317547305) node (311) {} -- (1.25,1.0825317547305) node (312) {} -- (0.9375,1.6237976320958) node (313) {}--cycle;\fill  (0.625,1.0825317547305) circle[radius=1.3pt]  (1.25,1.0825317547305)  circle[radius=1.3pt]  (0.9375,1.6237976320958) circle[radius=1.3pt] ;\draw  (1.25,1.0825317547305) node (321) {} -- (1.875,1.0825317547305) node (322) {} -- (1.5625,1.6237976320958) node (323) {}--cycle;\fill  (1.25,1.0825317547305) circle[radius=1.3pt]  (1.875,1.0825317547305)  circle[radius=1.3pt]  (1.5625,1.6237976320958) circle[radius=1.3pt] ;\draw  (0.9375,1.6237976320958) node (331) {} -- (1.5625,1.6237976320958) node (332) {} -- (1.25,2.1650635094611) node (333) {}--cycle;\fill  (0.9375,1.6237976320958) circle[radius=1.3pt]  (1.5625,1.6237976320958)  circle[radius=1.3pt]  (1.25,2.1650635094611) circle[radius=1.3pt] ;
\path  (111) node[below left] {2}  (112) node[below] {-1} (113) node[left] {-1}  (123) node[right] {1};
\end{tikzpicture}
\begin{tikzpicture}[scale=1.4]
\draw  (0,0) node (111) {} -- (0.625,0) node (112) {} -- (0.3125,0.54126587736527) node (113) {}--cycle;\fill  (0,0) circle[radius=1.3pt]  (0.625,0)  circle[radius=1.3pt]  (0.3125,0.54126587736527) circle[radius=1.3pt] ;\draw  (0.625,0) node (121) {} -- (1.25,0) node (122) {} -- (0.9375,0.54126587736527) node (123) {}--cycle;\fill  (0.625,0) circle[radius=1.3pt]  (1.25,0)  circle[radius=1.3pt]  (0.9375,0.54126587736527) circle[radius=1.3pt] ;\draw  (0.3125,0.54126587736527) node (131) {} -- (0.9375,0.54126587736527) node (132) {} -- (0.625,1.0825317547305) node (133) {}--cycle;\fill  (0.3125,0.54126587736527) circle[radius=1.3pt]  (0.9375,0.54126587736527)  circle[radius=1.3pt]  (0.625,1.0825317547305) circle[radius=1.3pt] ;\draw  (1.25,0) node (211) {} -- (1.875,0) node (212) {} -- (1.5625,0.54126587736527) node (213) {}--cycle;\fill  (1.25,0) circle[radius=1.3pt]  (1.875,0)  circle[radius=1.3pt]  (1.5625,0.54126587736527) circle[radius=1.3pt] ;\draw  (1.875,0) node (221) {} -- (2.5,0) node (222) {} -- (2.1875,0.54126587736527) node (223) {}--cycle;\fill  (1.875,0) circle[radius=1.3pt]  (2.5,0)  circle[radius=1.3pt]  (2.1875,0.54126587736527) circle[radius=1.3pt] ;\draw  (1.5625,0.54126587736527) node (231) {} -- (2.1875,0.54126587736527) node (232) {} -- (1.875,1.0825317547305) node (233) {}--cycle;\fill  (1.5625,0.54126587736527) circle[radius=1.3pt]  (2.1875,0.54126587736527)  circle[radius=1.3pt]  (1.875,1.0825317547305) circle[radius=1.3pt] ;\draw  (0.625,1.0825317547305) node (311) {} -- (1.25,1.0825317547305) node (312) {} -- (0.9375,1.6237976320958) node (313) {}--cycle;\fill  (0.625,1.0825317547305) circle[radius=1.3pt]  (1.25,1.0825317547305)  circle[radius=1.3pt]  (0.9375,1.6237976320958) circle[radius=1.3pt] ;\draw  (1.25,1.0825317547305) node (321) {} -- (1.875,1.0825317547305) node (322) {} -- (1.5625,1.6237976320958) node (323) {}--cycle;\fill  (1.25,1.0825317547305) circle[radius=1.3pt]  (1.875,1.0825317547305)  circle[radius=1.3pt]  (1.5625,1.6237976320958) circle[radius=1.3pt] ;\draw  (0.9375,1.6237976320958) node (331) {} -- (1.5625,1.6237976320958) node (332) {} -- (1.25,2.1650635094611) node (333) {}--cycle;\fill  (0.9375,1.6237976320958) circle[radius=1.3pt]  (1.5625,1.6237976320958)  circle[radius=1.3pt]  (1.25,2.1650635094611) circle[radius=1.3pt] ;
\path  (122)node[below] {2}  (121) node[below]{-1}  (212) node[below]{-1}  (123) node[above=2.5pt, right=-2pt]{-1}  (213) node[above=2.5pt, left=-2.5pt]{-1} (113)node[left]  {1}  (223) node[right]{1};
\end{tikzpicture}
\caption{\protect{Eigenfunctions of $V^G_1$ (left) and $V^G_2$ (center and right) with eigenvalue $-\frac32$.}}\label{fig:6series}
\end{figure}

A basis for this eigenspace can be chosen so as to be localized around points of  $V^G_{n-1}$. For $n\geq1$ it is generated by placing copies of the eigenfunction on $V^G_1$ shown on the left in Figure~\ref{fig:6series} on $n-1$ cells in $V^G_n$, so that the value $2$ occurs at a point of $x\in V^G_{n-1}$.  If $x\in V^G_0$ one copy suffices, as shown in the center of the figure. If $x\in V^G_{n-1}\setminus V^G_0$ two copies are needed, as shown on the right.  The functions are zero at all unlableled vertices in the figure.

Observe that the functions in Proposition~\ref{prop:Gnexcept1} exactly satisfy the condition given in Proposition~\ref{prop:Gndec} for the restriction of an eigenfunction $f$ of $\Delta^J_n$ with eigenvalue $-1$ to the set $V^G_n$.  Moreover, if $f|_{V^G_n}$ is such a function then $\Phi_3f$ vanishes on $V^H_n$, so coincides with $f_{V^H_n}$, also by Proposition~\ref{prop:Gndec}.  We have therefore extended the bijective correspondence of eigenspaces in Proposition~\ref{prop:Gndec} to the exceptional value. 

\begin{corollary}
$R_3:\sigma(\Delta^J_n)\to\sigma(\Delta^G_n)$ with a corresponding bijection of eigenspaces. From $J_n$ to $G_n$ the bijection  is simply restriction of an eigenfunction $f$ of $\Delta^J_n$ to $V^G_n$. From $G_n$ to $J_n$ the map is as follows. Suppose $g$ is an eigenfunction of $\Delta^G_n$ with eigenvalue $\zeta$. If $\zeta=-\frac32$ extend $g$ to $V^J_n$ so it vanishes on $V^H_n$. If $\zeta\neq-\frac32$  and $z\in R_3^{-1}(\zeta)$ is one of the corresponding elements of $\sigma(\Delta^J_n)\setminus\{-1\}$ then define $f=g$ on $V^G_n$ and $f=(1+z)^{-1}\Phi_3g$ on $V^H_n$.
\end{corollary}

\begin{corollary}\label{cor:kerPhi3}
The kernel of $\Phi_3$ is the eigenspace of $\Delta^G_n$ with eigenvalue $-\frac32$.
\end{corollary}

Since our method for obtaining eigenfunctions of $\Delta^H_n$ from those of $\Delta^G_n$ (in Theorem~\ref{thm:main}) is to apply $\Phi_3$, Corollary~\ref{cor:kerPhi3} makes it clear that nothing of the kind is possible for the $-\frac32$ eigenspace of $\Delta^G_n$. However, we can use it to get eigenfunctions of $\Delta^H_n$ with eigenvalue $-1$.

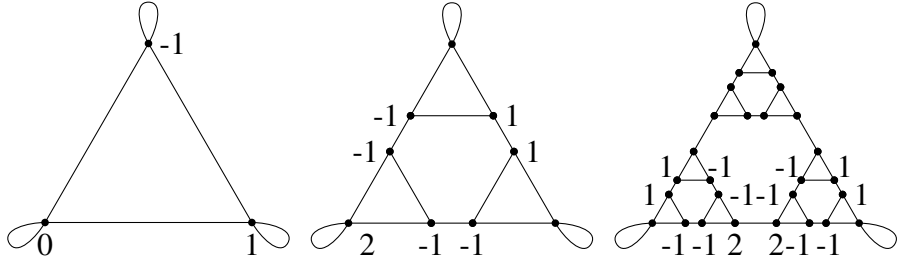
\begin{figure}
\begin{tikzpicture}[scale=1.1]
\draw  (0,0) node (1) {} -- (2.5,0) node (2) {} -- (1.25,2.1650635094611) node (3) {}--cycle;\fill  (0,0) circle[radius=1.3pt]  (2.5,0)  circle[radius=1.3pt]  (1.25,2.1650635094611) circle[radius=1.3pt] ;
\draw (0,0) to  [out=180,in=120]  (-.433,-.25)  to [out=300, in=240] (0,0);
\draw (2.5,0) to [out=0,in=60]   (2.933,-.25)  to [out=240, in=300] (2.5,0);
\draw (1.25,2.16506) to  [out=60,in=0]  (1.25,2.66506)  to  [out=180, in=120] (1.25,2.16506);
\path (1) node[below] {0} (2) node[below]{1} (3) node[right]{-1};
\end{tikzpicture}
\begin{tikzpicture}[scale=1.1]
\draw (1,0) node{} -- (1.5,0);\draw (0.5,0.86602540378444) node{} -- (0.75,1.2990381056767);\draw (2,0.86602540378444) node{} -- (1.75,1.2990381056767);\fill  (1,0) circle[radius=1.3pt]  (1.5,0)  circle[radius=1.3pt]  (0.5,0.86602540378444) circle[radius=1.3pt]  (0.75,1.2990381056767) circle[radius=1.3pt]   (2,0.86602540378444) circle[radius=1.3pt]  (1.75,1.2990381056767) circle[radius=1.3pt]     ;\draw  (0,0) node (11) {} -- (1,0) node (12) {} -- (0.5,0.86602540378444) node (13) {}--cycle;\fill  (0,0) circle[radius=1.3pt]  (1,0)  circle[radius=1.3pt]  (0.5,0.86602540378444) circle[radius=1.3pt] ;\draw  (1.5,0) node (21) {} -- (2.5,0) node (22) {} -- (2,0.86602540378444) node (23) {}--cycle;\fill  (1.5,0) circle[radius=1.3pt]  (2.5,0)  circle[radius=1.3pt]  (2,0.86602540378444) circle[radius=1.3pt] ;\draw  (0.75,1.2990381056767) node (31) {} -- (1.75,1.2990381056767) node (32) {} -- (1.25,2.1650635094611) node (33) {}--cycle;\fill  (0.75,1.2990381056767) circle[radius=1.3pt]  (1.75,1.2990381056767)  circle[radius=1.3pt]  (1.25,2.1650635094611) circle[radius=1.3pt] ;
\draw (0,0) to  [out=180,in=120]  (-.433,-.25)  to [out=300, in=240] (0,0);
\draw (2.5,0) to [out=0,in=60]   (2.933,-.25)  to [out=240, in=300] (2.5,0);
\draw (1.25,2.16506) to  [out=60,in=0]  (1.25,2.66506)  to  [out=180, in=120] (1.25,2.16506);
\path (11) node[below right]{2} (12) node[below]{-1} (21) node[below]{-1} (23) node[right]{1} (32) node[right] {1} (31) node[left]{-1} (13) node[left]{-1};
\end{tikzpicture}
\begin{tikzpicture}[scale=1.1]
\draw (1,0) node{} -- (1.5,0);\draw (0.5,0.86602540378444) node{} -- (0.75,1.2990381056767);\draw (2,0.86602540378444) node{} -- (1.75,1.2990381056767);\fill  (1,0) circle[radius=1.3pt]  (1.5,0)  circle[radius=1.3pt]  (0.5,0.86602540378444) circle[radius=1.3pt]  (0.75,1.2990381056767) circle[radius=1.3pt]   (2,0.86602540378444) circle[radius=1.3pt]  (1.75,1.2990381056767) circle[radius=1.3pt]     ;\draw (0.4,0) node{} -- (0.6,0);\draw (0.2,0.34641016151378) node{} -- (0.3,0.51961524227066);\draw (0.8,0.34641016151378) node{} -- (0.7,0.51961524227066);\fill  (0.4,0) circle[radius=1.3pt]  (0.6,0)  circle[radius=1.3pt]  (0.2,0.34641016151378) circle[radius=1.3pt]  (0.3,0.51961524227066) circle[radius=1.3pt]   (0.8,0.34641016151378) circle[radius=1.3pt]  (0.7,0.51961524227066) circle[radius=1.3pt]     ;\draw  (0,0) node (111) {} -- (0.4,0) node (112) {} -- (0.2,0.34641016151378) node (113) {}--cycle;\fill  (0,0) circle[radius=1.3pt]  (0.4,0)  circle[radius=1.3pt]  (0.2,0.34641016151378) circle[radius=1.3pt] ;\draw  (0.6,0) node (121) {} -- (1,0) node (122) {} -- (0.8,0.34641016151378) node (123) {}--cycle;\fill  (0.6,0) circle[radius=1.3pt]  (1,0)  circle[radius=1.3pt]  (0.8,0.34641016151378) circle[radius=1.3pt] ;\draw  (0.3,0.51961524227066) node (131) {} -- (0.7,0.51961524227066) node (132) {} -- (0.5,0.86602540378444) node (133) {}--cycle;\fill  (0.3,0.51961524227066) circle[radius=1.3pt]  (0.7,0.51961524227066)  circle[radius=1.3pt]  (0.5,0.86602540378444) circle[radius=1.3pt] ;\draw (1.9,0) node{} -- (2.1,0);\draw (1.7,0.34641016151378) node{} -- (1.8,0.51961524227066);\draw (2.3,0.34641016151378) node{} -- (2.2,0.51961524227066);\fill  (1.9,0) circle[radius=1.3pt]  (2.1,0)  circle[radius=1.3pt]  (1.7,0.34641016151378) circle[radius=1.3pt]  (1.8,0.51961524227066) circle[radius=1.3pt]   (2.3,0.34641016151378) circle[radius=1.3pt]  (2.2,0.51961524227066) circle[radius=1.3pt]     ;\draw  (1.5,0) node (211) {} -- (1.9,0) node (212) {} -- (1.7,0.34641016151378) node (213) {}--cycle;\fill  (1.5,0) circle[radius=1.3pt]  (1.9,0)  circle[radius=1.3pt]  (1.7,0.34641016151378) circle[radius=1.3pt] ;\draw  (2.1,0) node (221) {} -- (2.5,0) node (222) {} -- (2.3,0.34641016151378) node (223) {}--cycle;\fill  (2.1,0) circle[radius=1.3pt]  (2.5,0)  circle[radius=1.3pt]  (2.3,0.34641016151378) circle[radius=1.3pt] ;\draw  (1.8,0.51961524227066) node (231) {} -- (2.2,0.51961524227066) node (232) {} -- (2,0.86602540378444) node (233) {}--cycle;\fill  (1.8,0.51961524227066) circle[radius=1.3pt]  (2.2,0.51961524227066)  circle[radius=1.3pt]  (2,0.86602540378444) circle[radius=1.3pt] ;\draw (1.15,1.2990381056767) node{} -- (1.35,1.2990381056767);\draw (0.95,1.6454482671904) node{} -- (1.05,1.8186533479473);\draw (1.55,1.6454482671904) node{} -- (1.45,1.8186533479473);\fill  (1.15,1.2990381056767) circle[radius=1.3pt]  (1.35,1.2990381056767)  circle[radius=1.3pt]  (0.95,1.6454482671904) circle[radius=1.3pt]  (1.05,1.8186533479473) circle[radius=1.3pt]   (1.55,1.6454482671904) circle[radius=1.3pt]  (1.45,1.8186533479473) circle[radius=1.3pt]     ;\draw  (0.75,1.2990381056767) node (311) {} -- (1.15,1.2990381056767) node (312) {} -- (0.95,1.6454482671904) node (313) {}--cycle;\fill  (0.75,1.2990381056767) circle[radius=1.3pt]  (1.15,1.2990381056767)  circle[radius=1.3pt]  (0.95,1.6454482671904) circle[radius=1.3pt] ;\draw  (1.35,1.2990381056767) node (321) {} -- (1.75,1.2990381056767) node (322) {} -- (1.55,1.6454482671904) node (323) {}--cycle;\fill  (1.35,1.2990381056767) circle[radius=1.3pt]  (1.75,1.2990381056767)  circle[radius=1.3pt]  (1.55,1.6454482671904) circle[radius=1.3pt] ;\draw  (1.05,1.8186533479473) node (331) {} -- (1.45,1.8186533479473) node (332) {} -- (1.25,2.1650635094611) node (333) {}--cycle;\fill  (1.05,1.8186533479473) circle[radius=1.3pt]  (1.45,1.8186533479473)  circle[radius=1.3pt]  (1.25,2.1650635094611) circle[radius=1.3pt] ;
\draw (0,0) to  [out=180,in=120]  (-.433,-.25)  to [out=300, in=240] (0,0);
\draw (2.5,0) to [out=0,in=60]   (2.933,-.25)  to [out=240, in=300] (2.5,0);
\draw (1.25,2.16506) to  [out=60,in=0]  (1.25,2.66506)  to  [out=180, in=120] (1.25,2.16506);
\path (112) node[below=8pt, left=-5pt]{-1} (121) node[below=8pt, right=-8pt]{-1} (122) node[below]{2}  (211) node[below]{2} (212) node[below=8pt,left=-5pt]{-1} (221) node[below=8pt,right=-8pt]{-1}  (123) node[right]{-1} (132) node[right=4pt, above=-3pt] {-1} (131) node[left=4pt, above=-3pt]{1} (113) node[left]{1} (213) node[left]{-1} (231) node[left=5pt, above=-3pt] {-1} (223) node[right] {1} (232) node[right=4pt, above=-3pt]{1}  ;
\end{tikzpicture}
\caption{Eigenfunctions of $\Delta^H_n$ with eigenvalue $-1$ for $n=1,2,3$}\label{fig:Hn-1efn}
\end{figure}

\begin{proposition}\label{prop:Hn-1}
For $n\geq1$ the exceptional value $-1$ of the decimation in Theorem~\ref{thm:Hndecimation} is an eigenvalue of $\Delta^H_n$ with multiplicity  $\frac12(3^{n-1}+3)$.  If $n=1$ a basis consists of the function on the left in Figure~\ref{fig:Hn-1efn} and one rotate thereof. For $n\geq2$ a basis can be obtained from copies of the function on $H_2$ shown at center in Figure~\ref{fig:Hn-1efn} scaled to cells of scale $n-1$ and arranged as follows. A single copy may be placed so the value $2$ occurs at the attachment of a loop, or two copies may be placed on adjacent cells such that both endpoints of the connecting edge  carry the value $2$, as shown at right in Figure~\ref{fig:Hn-1efn}.
\end{proposition}

\begin{proof}
By Proposition~\ref{prop:Gnexcept1} we know $-\frac32$ is an eigenvalue of $\Delta^G_{n-1}$ with multiplicity $\frac12(3^{n-1}+3)$. Now $R_G(-\frac34)=-1$, so by Theorem~\ref{thm:SGspect} we have that $-\frac34$ is an eigenvalue of $\Delta^G_n$, and there is a bijection between the eigenspaces.  By Theorem~\ref{thm:main} then $-1$ is an eigenvalue of $\Delta^H_n$; again there is a bijection between the eigenspaces from which we deduce the stated multiplicity. 

One can obtain the eigenfunctions by taking values on $V^G_{n-1}$ as illustrated in Figure~\ref{fig:6series} and applying~\eqref{eq:Gefnextension}, calling the result $f$. Recall from Proposition~\ref{prop:Gnexcept1} that then $\sum_{j=1}^3 f(F_w(p_j))=0$ for all words $w$ of length $n$.  We compute from this constraint on the values and the requirement that $z=-\frac34$ that the value at a point in $V^G_n\setminus V^G_{n-1}$ is the average of the values at its neighbors. 

When $n=1$ we can reason as above for the eigenfunction of $\Delta^G_0$ having values $1,-1,0$ on $V_0$. Multiplying by $2$
gives the function on the left in Figure~\ref{fig:Hn-1efn}, and any two of the three rotates are linearly independent. 
When $n=2$ we do the same for the left function in Figure~\ref{fig:6series}. Then the above reasoning followed by the map $2\Phi_3$ produces the function on $H_2$ shown at the center in Figure~\ref{fig:Hn-1efn}.   For any $n\geq2$ copies of this can be rescaled to cells of size $n-1$ and placed so that  the common value $2$ occurs at the attachment point of a loop. If $n\geq3$ it is also possible to place two copies so the common value $2$ occurs at both ends of one of the $\frac12(3^{n-1}-3)$ edges in $H_n$ that were already present in $H_{n-1}$, see at right in Figure~\ref{fig:Hn-1efn}. The fact that these are linearly independent follows from the previously stated bijectivity or can readily be verified by hand.
\end{proof}

\begin{remark}
An easy consequence of the preceding is that the values of the $-1$ eigenfunctions of $\Delta^H_n$ sum to zero over all triangles of scale $n$ in the graph, though this is not the only constraint on such eigenfunctions.
\end{remark}

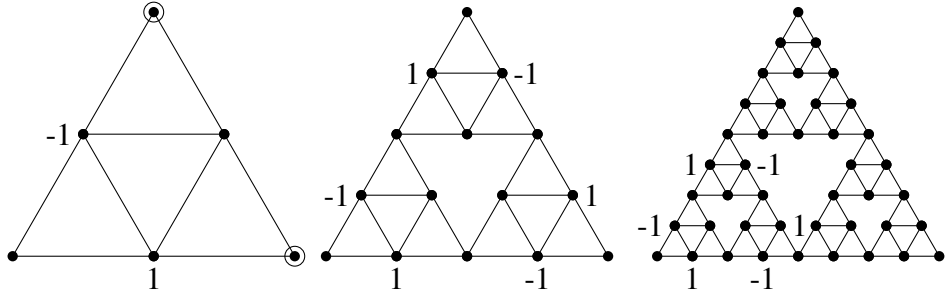
\begin{figure}
\begin{tikzpicture}[scale=1.5]
\draw  (0,0) node (11) {} -- (1.25,0) node (12) {} -- (0.625,1.0825317547305) node (13) {}--cycle;\fill  (0,0) circle[radius=1.3pt]  (1.25,0)  circle[radius=1.3pt]  (0.625,1.0825317547305) circle[radius=1.3pt] ;\draw  (1.25,0) node (21) {} -- (2.5,0) node (22) {} -- (1.875,1.0825317547305) node (23) {}--cycle;\fill  (1.25,0) circle[radius=1.3pt]  (2.5,0)  circle[radius=1.3pt]  (1.875,1.0825317547305) circle[radius=1.3pt] ;\draw  (0.625,1.0825317547305) node (31) {} -- (1.875,1.0825317547305) node (32) {} -- (1.25,2.1650635094611) node (33) {}--cycle;\fill  (0.625,1.0825317547305) circle[radius=1.3pt]  (1.875,1.0825317547305)  circle[radius=1.3pt]  (1.25,2.1650635094611) circle[radius=1.3pt] ;
\path (13) node[left]{-1} (12) node[below]{1};
\draw (33) circle[radius=2.5pt]  (22) circle[radius=2.5pt];
\end{tikzpicture}
\begin{tikzpicture}[scale=1.5]
\draw  (0,0) node (111) {} -- (0.625,0) node (112) {} -- (0.3125,0.54126587736527) node (113) {}--cycle;\fill  (0,0) circle[radius=1.3pt]  (0.625,0)  circle[radius=1.3pt]  (0.3125,0.54126587736527) circle[radius=1.3pt] ;\draw  (0.625,0) node (121) {} -- (1.25,0) node (122) {} -- (0.9375,0.54126587736527) node (123) {}--cycle;\fill  (0.625,0) circle[radius=1.3pt]  (1.25,0)  circle[radius=1.3pt]  (0.9375,0.54126587736527) circle[radius=1.3pt] ;\draw  (0.3125,0.54126587736527) node (131) {} -- (0.9375,0.54126587736527) node (132) {} -- (0.625,1.0825317547305) node (133) {}--cycle;\fill  (0.3125,0.54126587736527) circle[radius=1.3pt]  (0.9375,0.54126587736527)  circle[radius=1.3pt]  (0.625,1.0825317547305) circle[radius=1.3pt] ;\draw  (1.25,0) node (211) {} -- (1.875,0) node (212) {} -- (1.5625,0.54126587736527) node (213) {}--cycle;\fill  (1.25,0) circle[radius=1.3pt]  (1.875,0)  circle[radius=1.3pt]  (1.5625,0.54126587736527) circle[radius=1.3pt] ;\draw  (1.875,0) node (221) {} -- (2.5,0) node (222) {} -- (2.1875,0.54126587736527) node (223) {}--cycle;\fill  (1.875,0) circle[radius=1.3pt]  (2.5,0)  circle[radius=1.3pt]  (2.1875,0.54126587736527) circle[radius=1.3pt] ;\draw  (1.5625,0.54126587736527) node (231) {} -- (2.1875,0.54126587736527) node (232) {} -- (1.875,1.0825317547305) node (233) {}--cycle;\fill  (1.5625,0.54126587736527) circle[radius=1.3pt]  (2.1875,0.54126587736527)  circle[radius=1.3pt]  (1.875,1.0825317547305) circle[radius=1.3pt] ;\draw  (0.625,1.0825317547305) node (311) {} -- (1.25,1.0825317547305) node (312) {} -- (0.9375,1.6237976320958) node (313) {}--cycle;\fill  (0.625,1.0825317547305) circle[radius=1.3pt]  (1.25,1.0825317547305)  circle[radius=1.3pt]  (0.9375,1.6237976320958) circle[radius=1.3pt] ;\draw  (1.25,1.0825317547305) node (321) {} -- (1.875,1.0825317547305) node (322) {} -- (1.5625,1.6237976320958) node (323) {}--cycle;\fill  (1.25,1.0825317547305) circle[radius=1.3pt]  (1.875,1.0825317547305)  circle[radius=1.3pt]  (1.5625,1.6237976320958) circle[radius=1.3pt] ;\draw  (0.9375,1.6237976320958) node (331) {} -- (1.5625,1.6237976320958) node (332) {} -- (1.25,2.1650635094611) node (333) {}--cycle;\fill  (0.9375,1.6237976320958) circle[radius=1.3pt]  (1.5625,1.6237976320958)  circle[radius=1.3pt]  (1.25,2.1650635094611) circle[radius=1.3pt] ;
\path (112) node[below]{1} (212) node[below] {-1} (223) node[right]{1} (323) node[right]{-1} (313) node[left]{1} (113) node[left]{-1};
\end{tikzpicture}
\begin{tikzpicture}[scale=1.5]
\draw  (0,0) node (1111) {} -- (0.3125,0) node (1112) {} -- (0.15625,0.27063293868264) node (1113) {}--cycle;\fill  (0,0) circle[radius=1.3pt]  (0.3125,0)  circle[radius=1.3pt]  (0.15625,0.27063293868264) circle[radius=1.3pt] ;\draw  (0.3125,0) node (1121) {} -- (0.625,0) node (1122) {} -- (0.46875,0.27063293868264) node (1123) {}--cycle;\fill  (0.3125,0) circle[radius=1.3pt]  (0.625,0)  circle[radius=1.3pt]  (0.46875,0.27063293868264) circle[radius=1.3pt] ;\draw  (0.15625,0.27063293868264) node (1131) {} -- (0.46875,0.27063293868264) node (1132) {} -- (0.3125,0.54126587736527) node (1133) {}--cycle;\fill  (0.15625,0.27063293868264) circle[radius=1.3pt]  (0.46875,0.27063293868264)  circle[radius=1.3pt]  (0.3125,0.54126587736527) circle[radius=1.3pt] ;\draw  (0.625,0) node (1211) {} -- (0.9375,0) node (1212) {} -- (0.78125,0.27063293868264) node (1213) {}--cycle;\fill  (0.625,0) circle[radius=1.3pt]  (0.9375,0)  circle[radius=1.3pt]  (0.78125,0.27063293868264) circle[radius=1.3pt] ;\draw  (0.9375,0) node (1221) {} -- (1.25,0) node (1222) {} -- (1.09375,0.27063293868264) node (1223) {}--cycle;\fill  (0.9375,0) circle[radius=1.3pt]  (1.25,0)  circle[radius=1.3pt]  (1.09375,0.27063293868264) circle[radius=1.3pt] ;\draw  (0.78125,0.27063293868264) node (1231) {} -- (1.09375,0.27063293868264) node (1232) {} -- (0.9375,0.54126587736527) node (1233) {}--cycle;\fill  (0.78125,0.27063293868264) circle[radius=1.3pt]  (1.09375,0.27063293868264)  circle[radius=1.3pt]  (0.9375,0.54126587736527) circle[radius=1.3pt] ;\draw  (0.3125,0.54126587736527) node (1311) {} -- (0.625,0.54126587736527) node (1312) {} -- (0.46875,0.81189881604791) node (1313) {}--cycle;\fill  (0.3125,0.54126587736527) circle[radius=1.3pt]  (0.625,0.54126587736527)  circle[radius=1.3pt]  (0.46875,0.81189881604791) circle[radius=1.3pt] ;\draw  (0.625,0.54126587736527) node (1321) {} -- (0.9375,0.54126587736527) node (1322) {} -- (0.78125,0.81189881604791) node (1323) {}--cycle;\fill  (0.625,0.54126587736527) circle[radius=1.3pt]  (0.9375,0.54126587736527)  circle[radius=1.3pt]  (0.78125,0.81189881604791) circle[radius=1.3pt] ;\draw  (0.46875,0.81189881604791) node (1331) {} -- (0.78125,0.81189881604791) node (1332) {} -- (0.625,1.0825317547305) node (1333) {}--cycle;\fill  (0.46875,0.81189881604791) circle[radius=1.3pt]  (0.78125,0.81189881604791)  circle[radius=1.3pt]  (0.625,1.0825317547305) circle[radius=1.3pt] ;\draw  (1.25,0) node (2111) {} -- (1.5625,0) node (2112) {} -- (1.40625,0.27063293868264) node (2113) {}--cycle;\fill  (1.25,0) circle[radius=1.3pt]  (1.5625,0)  circle[radius=1.3pt]  (1.40625,0.27063293868264) circle[radius=1.3pt] ;\draw  (1.5625,0) node (2121) {} -- (1.875,0) node (2122) {} -- (1.71875,0.27063293868264) node (2123) {}--cycle;\fill  (1.5625,0) circle[radius=1.3pt]  (1.875,0)  circle[radius=1.3pt]  (1.71875,0.27063293868264) circle[radius=1.3pt] ;\draw  (1.40625,0.27063293868264) node (2131) {} -- (1.71875,0.27063293868264) node (2132) {} -- (1.5625,0.54126587736527) node (2133) {}--cycle;\fill  (1.40625,0.27063293868264) circle[radius=1.3pt]  (1.71875,0.27063293868264)  circle[radius=1.3pt]  (1.5625,0.54126587736527) circle[radius=1.3pt] ;\draw  (1.875,0) node (2211) {} -- (2.1875,0) node (2212) {} -- (2.03125,0.27063293868264) node (2213) {}--cycle;\fill  (1.875,0) circle[radius=1.3pt]  (2.1875,0)  circle[radius=1.3pt]  (2.03125,0.27063293868264) circle[radius=1.3pt] ;\draw  (2.1875,0) node (2221) {} -- (2.5,0) node (2222) {} -- (2.34375,0.27063293868264) node (2223) {}--cycle;\fill  (2.1875,0) circle[radius=1.3pt]  (2.5,0)  circle[radius=1.3pt]  (2.34375,0.27063293868264) circle[radius=1.3pt] ;\draw  (2.03125,0.27063293868264) node (2231) {} -- (2.34375,0.27063293868264) node (2232) {} -- (2.1875,0.54126587736527) node (2233) {}--cycle;\fill  (2.03125,0.27063293868264) circle[radius=1.3pt]  (2.34375,0.27063293868264)  circle[radius=1.3pt]  (2.1875,0.54126587736527) circle[radius=1.3pt] ;\draw  (1.5625,0.54126587736527) node (2311) {} -- (1.875,0.54126587736527) node (2312) {} -- (1.71875,0.81189881604791) node (2313) {}--cycle;\fill  (1.5625,0.54126587736527) circle[radius=1.3pt]  (1.875,0.54126587736527)  circle[radius=1.3pt]  (1.71875,0.81189881604791) circle[radius=1.3pt] ;\draw  (1.875,0.54126587736527) node (2321) {} -- (2.1875,0.54126587736527) node (2322) {} -- (2.03125,0.81189881604791) node (2323) {}--cycle;\fill  (1.875,0.54126587736527) circle[radius=1.3pt]  (2.1875,0.54126587736527)  circle[radius=1.3pt]  (2.03125,0.81189881604791) circle[radius=1.3pt] ;\draw  (1.71875,0.81189881604791) node (2331) {} -- (2.03125,0.81189881604791) node (2332) {} -- (1.875,1.0825317547305) node (2333) {}--cycle;\fill  (1.71875,0.81189881604791) circle[radius=1.3pt]  (2.03125,0.81189881604791)  circle[radius=1.3pt]  (1.875,1.0825317547305) circle[radius=1.3pt] ;\draw  (0.625,1.0825317547305) node (3111) {} -- (0.9375,1.0825317547305) node (3112) {} -- (0.78125,1.3531646934132) node (3113) {}--cycle;\fill  (0.625,1.0825317547305) circle[radius=1.3pt]  (0.9375,1.0825317547305)  circle[radius=1.3pt]  (0.78125,1.3531646934132) circle[radius=1.3pt] ;\draw  (0.9375,1.0825317547305) node (3121) {} -- (1.25,1.0825317547305) node (3122) {} -- (1.09375,1.3531646934132) node (3123) {}--cycle;\fill  (0.9375,1.0825317547305) circle[radius=1.3pt]  (1.25,1.0825317547305)  circle[radius=1.3pt]  (1.09375,1.3531646934132) circle[radius=1.3pt] ;\draw  (0.78125,1.3531646934132) node (3131) {} -- (1.09375,1.3531646934132) node (3132) {} -- (0.9375,1.6237976320958) node (3133) {}--cycle;\fill  (0.78125,1.3531646934132) circle[radius=1.3pt]  (1.09375,1.3531646934132)  circle[radius=1.3pt]  (0.9375,1.6237976320958) circle[radius=1.3pt] ;\draw  (1.25,1.0825317547305) node (3211) {} -- (1.5625,1.0825317547305) node (3212) {} -- (1.40625,1.3531646934132) node (3213) {}--cycle;\fill  (1.25,1.0825317547305) circle[radius=1.3pt]  (1.5625,1.0825317547305)  circle[radius=1.3pt]  (1.40625,1.3531646934132) circle[radius=1.3pt] ;\draw  (1.5625,1.0825317547305) node (3221) {} -- (1.875,1.0825317547305) node (3222) {} -- (1.71875,1.3531646934132) node (3223) {}--cycle;\fill  (1.5625,1.0825317547305) circle[radius=1.3pt]  (1.875,1.0825317547305)  circle[radius=1.3pt]  (1.71875,1.3531646934132) circle[radius=1.3pt] ;\draw  (1.40625,1.3531646934132) node (3231) {} -- (1.71875,1.3531646934132) node (3232) {} -- (1.5625,1.6237976320958) node (3233) {}--cycle;\fill  (1.40625,1.3531646934132) circle[radius=1.3pt]  (1.71875,1.3531646934132)  circle[radius=1.3pt]  (1.5625,1.6237976320958) circle[radius=1.3pt] ;\draw  (0.9375,1.6237976320958) node (3311) {} -- (1.25,1.6237976320958) node (3312) {} -- (1.09375,1.8944305707785) node (3313) {}--cycle;\fill  (0.9375,1.6237976320958) circle[radius=1.3pt]  (1.25,1.6237976320958)  circle[radius=1.3pt]  (1.09375,1.8944305707785) circle[radius=1.3pt] ;\draw  (1.25,1.6237976320958) node (3321) {} -- (1.5625,1.6237976320958) node (3322) {} -- (1.40625,1.8944305707785) node (3323) {}--cycle;\fill  (1.25,1.6237976320958) circle[radius=1.3pt]  (1.5625,1.6237976320958)  circle[radius=1.3pt]  (1.40625,1.8944305707785) circle[radius=1.3pt] ;\draw  (1.09375,1.8944305707785) node (3331) {} -- (1.40625,1.8944305707785) node (3332) {} -- (1.25,2.1650635094611) node (3333) {}--cycle;\fill  (1.09375,1.8944305707785) circle[radius=1.3pt]  (1.40625,1.8944305707785)  circle[radius=1.3pt]  (1.25,2.1650635094611) circle[radius=1.3pt] ;
\path (1112) node[below]{1} (1212) node[below]{-1} (1223) node[right]{1} (1323) node[right]{-1} (1313) node[left]{1} (1113) node[left]{-1};
\end{tikzpicture}
\caption{Eigenfunctions of $\Delta^G_n$ with eigenvalue $-\frac54$ for $n=2$ (center) and $n=3$ (right) are constructed from the function at left, which satisfies the eigenfunction equation at all but the circled points, by chaining around holes in the graph.}\label{fig:5series}
\end{figure}

The eigenfunctions corresponding to the exceptional value $-\frac54$ for the decimation of $\Delta^G_n$ are known to be in one-to-one correspondence with certain loops on the cell graph. The left image in Figure~\ref{fig:5series} shows a function on $V^G_1$ by listing the values at vertices, with the vertices where the function value is zero left unmarked.  This function satisfies the eigenfunction equation with eigenvalue $-\frac54$ at all vertices except the two that are circled. By symmetry, one may join copies of this function  along non-trivial cycles so as to obtain eigenfunctions of $\Delta^G_n$ with eigenvalue $-\frac54$ for any $n\geq 2$. The cycles surround ``holes'' in the cell graph with size at least one scale larger than a cell, as illustrated on the center and right in Figure~\ref{fig:5series}. Counting the holes gives the multiplicity of the eigenspace.

\begin{proposition}[\protect{\cite{DSV}, see also~\cite[Chapter 3]{Strichartzbook}}] \label{prop:Gnexcept2}
The exceptional value $-\frac54$ is an eigenvalue for $n\geq2$ with multiplicity $\frac12(3^{n-1}-1)$. An explicit basis for the eigenspace is obtained by chaining the  function on the left in Figure~\ref{fig:5series} along cycles around holes of scale at least $n-1$.
\end{proposition}

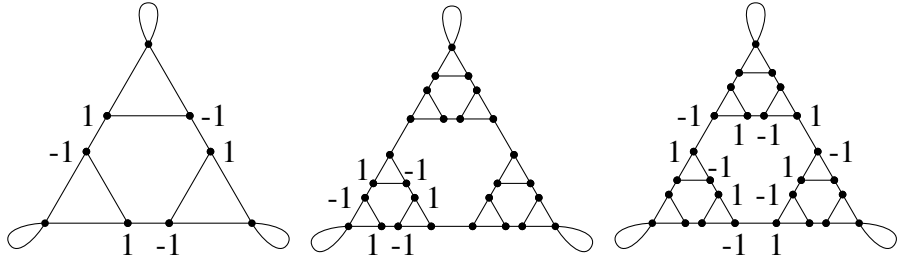
\begin{figure}
\begin{tikzpicture}[scale=1.1]
\draw (1,0) node{} -- (1.5,0);\draw (0.5,0.86602540378444) node{} -- (0.75,1.2990381056767);\draw (2,0.86602540378444) node{} -- (1.75,1.2990381056767);\fill  (1,0) circle[radius=1.3pt]  (1.5,0)  circle[radius=1.3pt]  (0.5,0.86602540378444) circle[radius=1.3pt]  (0.75,1.2990381056767) circle[radius=1.3pt]   (2,0.86602540378444) circle[radius=1.3pt]  (1.75,1.2990381056767) circle[radius=1.3pt]     ;\draw  (0,0) node (11) {} -- (1,0) node (12) {} -- (0.5,0.86602540378444) node (13) {}--cycle;\fill  (0,0) circle[radius=1.3pt]  (1,0)  circle[radius=1.3pt]  (0.5,0.86602540378444) circle[radius=1.3pt] ;\draw  (1.5,0) node (21) {} -- (2.5,0) node (22) {} -- (2,0.86602540378444) node (23) {}--cycle;\fill  (1.5,0) circle[radius=1.3pt]  (2.5,0)  circle[radius=1.3pt]  (2,0.86602540378444) circle[radius=1.3pt] ;\draw  (0.75,1.2990381056767) node (31) {} -- (1.75,1.2990381056767) node (32) {} -- (1.25,2.1650635094611) node (33) {}--cycle;\fill  (0.75,1.2990381056767) circle[radius=1.3pt]  (1.75,1.2990381056767)  circle[radius=1.3pt]  (1.25,2.1650635094611) circle[radius=1.3pt] ;
\draw (0,0) to  [out=180,in=120]  (-.433,-.25)  to [out=300, in=240] (0,0);
\draw (2.5,0) to [out=0,in=60]   (2.933,-.25)  to [out=240, in=300] (2.5,0);
\draw (1.25,2.16506) to  [out=60,in=0]  (1.25,2.66506)  to  [out=180, in=120] (1.25,2.16506);
\path (12) node[below]{1} (21) node[below]{-1} (23) node[right]{1} (32) node[right] {-1} (31) node[left]{1} (13) node[left]{-1};
\end{tikzpicture}
\begin{tikzpicture}[scale=1.1]
\draw (1,0) node{} -- (1.5,0);\draw (0.5,0.86602540378444) node{} -- (0.75,1.2990381056767);\draw (2,0.86602540378444) node{} -- (1.75,1.2990381056767);\fill  (1,0) circle[radius=1.3pt]  (1.5,0)  circle[radius=1.3pt]  (0.5,0.86602540378444) circle[radius=1.3pt]  (0.75,1.2990381056767) circle[radius=1.3pt]   (2,0.86602540378444) circle[radius=1.3pt]  (1.75,1.2990381056767) circle[radius=1.3pt]     ;\draw (0.4,0) node{} -- (0.6,0);\draw (0.2,0.34641016151378) node{} -- (0.3,0.51961524227066);\draw (0.8,0.34641016151378) node{} -- (0.7,0.51961524227066);\fill  (0.4,0) circle[radius=1.3pt]  (0.6,0)  circle[radius=1.3pt]  (0.2,0.34641016151378) circle[radius=1.3pt]  (0.3,0.51961524227066) circle[radius=1.3pt]   (0.8,0.34641016151378) circle[radius=1.3pt]  (0.7,0.51961524227066) circle[radius=1.3pt]     ;\draw  (0,0) node (111) {} -- (0.4,0) node (112) {} -- (0.2,0.34641016151378) node (113) {}--cycle;\fill  (0,0) circle[radius=1.3pt]  (0.4,0)  circle[radius=1.3pt]  (0.2,0.34641016151378) circle[radius=1.3pt] ;\draw  (0.6,0) node (121) {} -- (1,0) node (122) {} -- (0.8,0.34641016151378) node (123) {}--cycle;\fill  (0.6,0) circle[radius=1.3pt]  (1,0)  circle[radius=1.3pt]  (0.8,0.34641016151378) circle[radius=1.3pt] ;\draw  (0.3,0.51961524227066) node (131) {} -- (0.7,0.51961524227066) node (132) {} -- (0.5,0.86602540378444) node (133) {}--cycle;\fill  (0.3,0.51961524227066) circle[radius=1.3pt]  (0.7,0.51961524227066)  circle[radius=1.3pt]  (0.5,0.86602540378444) circle[radius=1.3pt] ;\draw (1.9,0) node{} -- (2.1,0);\draw (1.7,0.34641016151378) node{} -- (1.8,0.51961524227066);\draw (2.3,0.34641016151378) node{} -- (2.2,0.51961524227066);\fill  (1.9,0) circle[radius=1.3pt]  (2.1,0)  circle[radius=1.3pt]  (1.7,0.34641016151378) circle[radius=1.3pt]  (1.8,0.51961524227066) circle[radius=1.3pt]   (2.3,0.34641016151378) circle[radius=1.3pt]  (2.2,0.51961524227066) circle[radius=1.3pt]     ;\draw  (1.5,0) node (211) {} -- (1.9,0) node (212) {} -- (1.7,0.34641016151378) node (213) {}--cycle;\fill  (1.5,0) circle[radius=1.3pt]  (1.9,0)  circle[radius=1.3pt]  (1.7,0.34641016151378) circle[radius=1.3pt] ;\draw  (2.1,0) node (221) {} -- (2.5,0) node (222) {} -- (2.3,0.34641016151378) node (223) {}--cycle;\fill  (2.1,0) circle[radius=1.3pt]  (2.5,0)  circle[radius=1.3pt]  (2.3,0.34641016151378) circle[radius=1.3pt] ;\draw  (1.8,0.51961524227066) node (231) {} -- (2.2,0.51961524227066) node (232) {} -- (2,0.86602540378444) node (233) {}--cycle;\fill  (1.8,0.51961524227066) circle[radius=1.3pt]  (2.2,0.51961524227066)  circle[radius=1.3pt]  (2,0.86602540378444) circle[radius=1.3pt] ;\draw (1.15,1.2990381056767) node{} -- (1.35,1.2990381056767);\draw (0.95,1.6454482671904) node{} -- (1.05,1.8186533479473);\draw (1.55,1.6454482671904) node{} -- (1.45,1.8186533479473);\fill  (1.15,1.2990381056767) circle[radius=1.3pt]  (1.35,1.2990381056767)  circle[radius=1.3pt]  (0.95,1.6454482671904) circle[radius=1.3pt]  (1.05,1.8186533479473) circle[radius=1.3pt]   (1.55,1.6454482671904) circle[radius=1.3pt]  (1.45,1.8186533479473) circle[radius=1.3pt]     ;\draw  (0.75,1.2990381056767) node (311) {} -- (1.15,1.2990381056767) node (312) {} -- (0.95,1.6454482671904) node (313) {}--cycle;\fill  (0.75,1.2990381056767) circle[radius=1.3pt]  (1.15,1.2990381056767)  circle[radius=1.3pt]  (0.95,1.6454482671904) circle[radius=1.3pt] ;\draw  (1.35,1.2990381056767) node (321) {} -- (1.75,1.2990381056767) node (322) {} -- (1.55,1.6454482671904) node (323) {}--cycle;\fill  (1.35,1.2990381056767) circle[radius=1.3pt]  (1.75,1.2990381056767)  circle[radius=1.3pt]  (1.55,1.6454482671904) circle[radius=1.3pt] ;\draw  (1.05,1.8186533479473) node (331) {} -- (1.45,1.8186533479473) node (332) {} -- (1.25,2.1650635094611) node (333) {}--cycle;\fill  (1.05,1.8186533479473) circle[radius=1.3pt]  (1.45,1.8186533479473)  circle[radius=1.3pt]  (1.25,2.1650635094611) circle[radius=1.3pt] ;
\draw (0,0) to  [out=180,in=120]  (-.433,-.25)  to [out=300, in=240] (0,0);
\draw (2.5,0) to [out=0,in=60]   (2.933,-.25)  to [out=240, in=300] (2.5,0);
\draw (1.25,2.16506) to  [out=60,in=0]  (1.25,2.66506)  to  [out=180, in=120] (1.25,2.16506);
\path (112) node[below=7pt,left=-4pt]{1} (121) node[below=7pt, right=-7pt]{-1} (123) node[right]{1} (132) node[right=4pt, above=-3pt] {-1} (131) node[left=4pt, above=-3pt]{1} (113) node[left]{-1};
\end{tikzpicture}
\begin{tikzpicture}[scale=1.1]
\draw (1,0) node{} -- (1.5,0);\draw (0.5,0.86602540378444) node{} -- (0.75,1.2990381056767);\draw (2,0.86602540378444) node{} -- (1.75,1.2990381056767);\fill  (1,0) circle[radius=1.3pt]  (1.5,0)  circle[radius=1.3pt]  (0.5,0.86602540378444) circle[radius=1.3pt]  (0.75,1.2990381056767) circle[radius=1.3pt]   (2,0.86602540378444) circle[radius=1.3pt]  (1.75,1.2990381056767) circle[radius=1.3pt]     ;\draw (0.4,0) node{} -- (0.6,0);\draw (0.2,0.34641016151378) node{} -- (0.3,0.51961524227066);\draw (0.8,0.34641016151378) node{} -- (0.7,0.51961524227066);\fill  (0.4,0) circle[radius=1.3pt]  (0.6,0)  circle[radius=1.3pt]  (0.2,0.34641016151378) circle[radius=1.3pt]  (0.3,0.51961524227066) circle[radius=1.3pt]   (0.8,0.34641016151378) circle[radius=1.3pt]  (0.7,0.51961524227066) circle[radius=1.3pt]     ;\draw  (0,0) node (111) {} -- (0.4,0) node (112) {} -- (0.2,0.34641016151378) node (113) {}--cycle;\fill  (0,0) circle[radius=1.3pt]  (0.4,0)  circle[radius=1.3pt]  (0.2,0.34641016151378) circle[radius=1.3pt] ;\draw  (0.6,0) node (121) {} -- (1,0) node (122) {} -- (0.8,0.34641016151378) node (123) {}--cycle;\fill  (0.6,0) circle[radius=1.3pt]  (1,0)  circle[radius=1.3pt]  (0.8,0.34641016151378) circle[radius=1.3pt] ;\draw  (0.3,0.51961524227066) node (131) {} -- (0.7,0.51961524227066) node (132) {} -- (0.5,0.86602540378444) node (133) {}--cycle;\fill  (0.3,0.51961524227066) circle[radius=1.3pt]  (0.7,0.51961524227066)  circle[radius=1.3pt]  (0.5,0.86602540378444) circle[radius=1.3pt] ;\draw (1.9,0) node{} -- (2.1,0);\draw (1.7,0.34641016151378) node{} -- (1.8,0.51961524227066);\draw (2.3,0.34641016151378) node{} -- (2.2,0.51961524227066);\fill  (1.9,0) circle[radius=1.3pt]  (2.1,0)  circle[radius=1.3pt]  (1.7,0.34641016151378) circle[radius=1.3pt]  (1.8,0.51961524227066) circle[radius=1.3pt]   (2.3,0.34641016151378) circle[radius=1.3pt]  (2.2,0.51961524227066) circle[radius=1.3pt]     ;\draw  (1.5,0) node (211) {} -- (1.9,0) node (212) {} -- (1.7,0.34641016151378) node (213) {}--cycle;\fill  (1.5,0) circle[radius=1.3pt]  (1.9,0)  circle[radius=1.3pt]  (1.7,0.34641016151378) circle[radius=1.3pt] ;\draw  (2.1,0) node (221) {} -- (2.5,0) node (222) {} -- (2.3,0.34641016151378) node (223) {}--cycle;\fill  (2.1,0) circle[radius=1.3pt]  (2.5,0)  circle[radius=1.3pt]  (2.3,0.34641016151378) circle[radius=1.3pt] ;\draw  (1.8,0.51961524227066) node (231) {} -- (2.2,0.51961524227066) node (232) {} -- (2,0.86602540378444) node (233) {}--cycle;\fill  (1.8,0.51961524227066) circle[radius=1.3pt]  (2.2,0.51961524227066)  circle[radius=1.3pt]  (2,0.86602540378444) circle[radius=1.3pt] ;\draw (1.15,1.2990381056767) node{} -- (1.35,1.2990381056767);\draw (0.95,1.6454482671904) node{} -- (1.05,1.8186533479473);\draw (1.55,1.6454482671904) node{} -- (1.45,1.8186533479473);\fill  (1.15,1.2990381056767) circle[radius=1.3pt]  (1.35,1.2990381056767)  circle[radius=1.3pt]  (0.95,1.6454482671904) circle[radius=1.3pt]  (1.05,1.8186533479473) circle[radius=1.3pt]   (1.55,1.6454482671904) circle[radius=1.3pt]  (1.45,1.8186533479473) circle[radius=1.3pt]     ;\draw  (0.75,1.2990381056767) node (311) {} -- (1.15,1.2990381056767) node (312) {} -- (0.95,1.6454482671904) node (313) {}--cycle;\fill  (0.75,1.2990381056767) circle[radius=1.3pt]  (1.15,1.2990381056767)  circle[radius=1.3pt]  (0.95,1.6454482671904) circle[radius=1.3pt] ;\draw  (1.35,1.2990381056767) node (321) {} -- (1.75,1.2990381056767) node (322) {} -- (1.55,1.6454482671904) node (323) {}--cycle;\fill  (1.35,1.2990381056767) circle[radius=1.3pt]  (1.75,1.2990381056767)  circle[radius=1.3pt]  (1.55,1.6454482671904) circle[radius=1.3pt] ;\draw  (1.05,1.8186533479473) node (331) {} -- (1.45,1.8186533479473) node (332) {} -- (1.25,2.1650635094611) node (333) {}--cycle;\fill  (1.05,1.8186533479473) circle[radius=1.3pt]  (1.45,1.8186533479473)  circle[radius=1.3pt]  (1.25,2.1650635094611) circle[radius=1.3pt] ;
\draw (0,0) to  [out=180,in=120]  (-.433,-.25)  to [out=300, in=240] (0,0);
\draw (2.5,0) to [out=0,in=60]   (2.933,-.25)  to [out=240, in=300] (2.5,0);
\draw (1.25,2.16506) to  [out=60,in=0]  (1.25,2.66506)  to  [out=180, in=120] (1.25,2.16506);
\path  (122) node[below]{-1}  (123) node[right]{1} (132) node[right=4pt, above=-3pt] {-1} (133) node[left]{1} (311) node[left]{-1} (312) node[below=7pt,left=-4pt]{1} (321) node[below=7pt, right=-7pt]{-1} (322) node[right]{1} (233) node[right]{-1}  (231) node[left=5pt, above=-3pt] {1} (213) node[left]{-1} (211) node[below]{1};
\end{tikzpicture}
\caption{Eigenfunctions of $\Delta^H_n$ with eigenvalues $-\frac53$ for $n=2$ (left) and $n=3$ (center and right) form loops around holes  in the graph of scale at least $n-1$.} \label{fig:Hn5series}
\end{figure}

Figure~\ref{fig:Hn5series} shows some of these eigenfunctions. In particular, the functions at the left and center are obtained by applying $3\Phi_3$ to the functions at the center and right in  Figure~\ref{fig:5series}. Evidently each such function is supported on a cycle around a hole and has alternating values going around the cycle.  Applying the bijection from Theorem~\ref{thm:main} we obtain the following result.

\begin{corollary}\label{cor:Hn-54}
For $n\geq2$ the exceptional value $-\frac53\in E_H$ is an eigenvalue of $\Delta^H_n$ with multiplicity $\frac12(3^{n-1}-1)$. A basis for the eigenspace consists of alternating functions around cycles in $H_n$ of any size larger than the smallest cycles, so has the homology of $H_{n-1}$.
\end{corollary}

The following well-known theorem is included only  for completeness of the description of $\sigma(\Delta^G_n)$ and to emphasize the connection to Theorem~\ref{thm:estructure of Hn}. It is proved in the same manner as the latter by counting the eigenfunctions obtained in Propositions~\ref{prop:Gnexcept1} and~\ref{prop:Gnexcept2} and by decimation using Theorem~\ref{thm:SGspect}.

\begin{theorem}[\protect{\cite{RammalToulouse, FukushimaShima,DSV}, see also~\cite[Chapter 3]{Strichartzbook}}]\label{thm:estructure of Gn}
The spectrum of $\Delta^G_n$ is the set
\begin{equation*}
	\sigma(\Delta^G_n) = \{0\} \cup\biggl( \bigcup_{i=0}^n R_G^{-i}\Bigl(-\frac32\Bigr) \biggr) \cup \biggl( \bigcup_{j=0}^{n-2} R_G^{-j}\Bigl(-\frac54 \Bigr)\biggr)
\end{equation*}
The multiplicities and eigenspaces are as follows:
\begin{itemize}
\item $0$ has multiplicity $1$ with constant eigenfunction
\item  $R_G^{-i}(-\frac32)$ contains $2^{i-1}$ eigenvalues when $1\leq i\leq n$, each with multiplicity $\frac12(3^{n-i}+3)$, while $-\frac32$ has multiplicity $\frac12(3^n+3)$. A basis for the eigenspace of any eigenvalue in $R_G^{-i}(-\frac32)$ is obtained by taking the basis for the $-1$ eigenspace of $\Delta^G_{n-i}$ described in Proposition~\ref{prop:Gnexcept1} and extending to $G_n$ by $i$ iterations of~\eqref{eq:Gefnextension}.
\item $R_G^{-j}(-\frac54)$ contains $2^j$ eigenvalues when $0\leq j\leq n-2$, each with multiplicity $\frac12(3^{n-j-1}-1)$. A basis for the eigenspace is obtained by taking the basis for the $-\frac54$ eigenspace for $\Delta^G_{n-j}$ from Proposition~\ref{prop:Gnexcept2} and extending to $G_n$ by $j$ iterations of~\eqref{eq:Gefnextension}.
\end{itemize}
\end{theorem}

The corresponding theorem for $\Delta^H_n$ is our main result, stated previously as Theorem~\ref{thm:estructure of Hn}.
\begin{proof}[Proof of \protect{Theorem~\ref{thm:estructure of Hn}}]
For $n=0$ and $n=1$ it is easy to check the description is correct using Proposition~\ref{prop:Hn-1} and the fact that the number of vertices in $V^H_n$ is $3^n$. For $n\geq2$ we induct. 

Consider the eigenvalues and eigenfunctions  obtained by the decimation in Theorem~\ref{thm:Hndecimation}.  Assuming for the induction that the spectrum for $\Delta^H_{n-1}$ is as stated we see that the forward orbit of any $\zeta\in\Delta^H_{n-1}$ under $R_H$ hits the set $S=\{-\frac53,-1,0\}$ (indeed it hits the smaller set $\{-1,0\}$, but we do not need this fact).  We use this to check that $\zeta$ cannot be a critical value of $R_H$ and to understand when $R_H^{-1}(\zeta)$ can be in $E_H$.  A useful observation in this regard is that $z>0$ implies $R_H(z)>0$ and the orbit of $z$ cannot be in $S$.  Applying this to the critical value $-\frac{25}{12}$ of $R_H$ we see $R_H\Bigl(-\frac{25}{12}\Bigr)=\frac{125}{12}$, so $\zeta\in\Delta^H_{n-1}$ implies $\zeta$ is not the critical value and thus $R_H^{-1}(\zeta)$ contains exactly two points.  For the set $E_H$ we compute the forward orbits $R_H(-\frac53)=0$ and $R_H(0)=0$, so these fall onto the fixed point $0\in S$, while $R_H(-1)=R_H(-\frac23)=-2$ and also $R_H(-2)=2$, after which the orbit never hits $S$.

Our inductive hypothesis says that for $0\leq i\leq n-2$ there are $2^i$ eigenvalues of $\Delta^H_{n-1}$ of the form $R_H^{-i}(-1)$, each with multiplicity $\frac12(3^{n-2-i}+3)$. For each of these eigenvalues $\zeta$ the above reasoning ensures $R_H^{-1}(\zeta)$ contains exactly two points, neither of which is in $E_H$. Then Theorem~\ref{thm:Hndecimation} shows both of these are eigenvalues of $\Delta^H_n$ with the same multiplicity as $\zeta$ and eigenfunctions obtained by the decimation formula~\eqref{eq:Hnextension}. Accordingly, for $1\leq i\leq n-1$ there are $2^i$ eigenvalues of $\Delta^H_n$ of the form $R_H^{-i}(-1)$; the case $i=0$ is then obtained from Proposition~\ref{prop:Hn-1}. 

In the case $n=2$ we can count that the preceding gives $3^n-1$ eigenfunctions and the remaining one comes from  Corollary~\ref{cor:Hn-54}.  At the same time we have inductively that for $n\geq3$ and $0\leq j\leq n-3$ there are $2^j$ eigenvalues of $\Delta^H_{n-1}$ each with multiplicity $\frac12(3^{n-2-j}-1)$. For each of these eigenvalues $\zeta$ our earlier reasoning about the orbits of $E_H$ and the critical value ensures $R_H^{-1}(\zeta)$ contains exactly two points, neither of which is in $E_H$, and thus by Theorem~\ref{thm:Hndecimation} we obtain two eigenvalues of $\Delta^H_n$ of the form $R_H^{-1}(\zeta)$ with eigenfunctions from~\eqref{eq:Hnextension}. Thus for $1\leq j\leq n-2$ we have $2^j$ eigenvalues of $\Delta^H_n$ with multiplicity $\frac12(3^{n-1-j}-1)$, and  this is also true for $j=0$ by Corollary~\ref{cor:Hn-54}.

The remaining possible eigenvalue $\zeta\in\Delta^H_{n-1}$ is $\zeta=0$, with multiplicity $1$. In this case $R_H^{-1}(0)=\{-\frac53,0\}$ and Theorem~\ref{thm:Hndecimation} says only that this ensures $0$ is an eigenvalue of $\Delta^H_n$, also with multiplicity $1$. Evidently the corresponding eigenfunction is constant.

The preceding does not yet ensure we have a full description of the spectrum, as it is a priori possible that the other values in $E_H$ could be eigenvalues. This possibility is eliminated by counting. Supposing inductively that our description gives the whole spectrum for $\Delta^H_{n-1}$, observe that this consists of $3^{n-1}$ eigenvalues counted with multiplicity. All but the zero eigenvalue gave two new eigenvalues of $\Delta^H_n$ with preservation of multiplicity by spectral decimation; adding back in the zero eigenvalue we thus obtain $2(3^{n-1})-1$ eigenvalues. We also have eigenvalues $-1$ and $-\frac54$ from Proposition~\ref{prop:Hn-1} and Corollary~\ref{cor:Hn-54}, with (respectively) multiplicities $\frac12(3^{n-1}+3)$ and $\frac12(3^{n-1}-1)$. Together these sum to $3^n$, which is the number of vertices in $V^H_n$, so we have found all eigenvalues and eigenfunctions.

What remains is the topological interpretation of the eigenspaces, which simply corresponds to the fact that all eigenfunctions of the types we have described were constructed by iterated application of~\eqref{eq:Hnextension} to the functions described in Proposition~\ref{prop:Hn-1} and Corollary~\ref{cor:Hn-54}.
\end{proof}

\begin{remark}
In our proof of Theorem~\ref{thm:estructure of Hn} we used the fact that we knew  the eigenspaces and hence the multiplicity of eigenvalues of $\Delta^H_n$ for the exceptional values $-1$ and $-\frac53$ exactly. However it is apparent in the proof that we needed only lower bounds for these multiplicities: once these lower bounds match the upper bound obtained by counting the vertices the proof is complete.  This observation shows that once Theorem~\ref{thm:Hndecimation} is known it is perfectly possible to obtain Theorem~\ref{thm:estructure of Hn} as follows.  Construct the eigenfunctions of $\Delta^H_n$ by hand for small $n$ and observe the topological structures described in Proposition~\ref{prop:Hn-1} and Corollary~\ref{cor:Hn-54} that allow one to produce some eigenfunctions of $\Delta^H_n$ for any $n$ with these exceptional values. Without knowing whether this gives all eigenfunctions, one still obtains lower bounds on the dimensions of both the eigenspaces corresponding to the exceptional values $-1$ and $-\frac53$ and the eigenspaces obtained from Theorem~\ref{thm:Hndecimation}. The argument given in the proof of Theorem~\ref{thm:estructure of Hn} gives matching upper bounds, completing the proof.
\end{remark}

\section{\protect{The method of Grigorchuk and \v{S}uni\'{c}.}}\label{sec:GScompare}

We have already mentioned that Grigorchuk and \v{S}uni\'{c}~\cite{GrigorchukSunic} computed the spectra of the sequence $\Delta^H_n$ by a different method. Our work shares many features with theirs, so it seems useful to identify the most important difference. For this reason we sketch some aspects of their approach in the same language used for our results above.  Note that although they consider the adjacency matrix rather than the Laplacian we phrase our discussion in terms of the latter; there is no essential difference because the graph is $3$-regular.

Considering the sequence of graphs $H_n$ in Figure~\ref{fig:Hn} it is apparent that we could assign vertices $\tilde V^H_n$ so that $\tilde V^H_n\subset \tilde V^H_{n+1}$ and the graph $H_{n+1}$ is obtained by subdividing edges of $H_n$ at these new vertices and introducing new edges between certain pairs of the new vertices. From this perspective it seems natural to seek a spectral decimation from $\Delta^H_{n+1}$ to $\Delta^H_n$ by taking the Schur complement of the subspace of functions supported on $\tilde V^H_{n+1}\setminus \tilde V^H_{n}$.  Although they come to this point by considering the permutation representation of the Hanoi towers group on the $n$-th level of a rooted $3$-tree, this Schur complement is the first step in the spectral analysis in~\cite{GrigorchukSunic}.

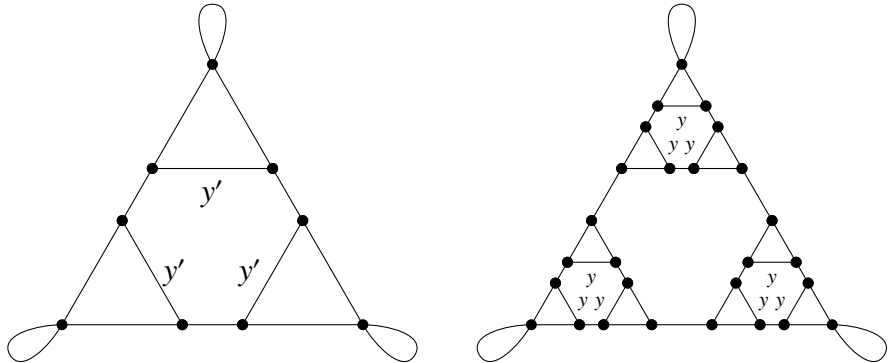
\begin{figure}
\begin{tikzpicture}[scale=1.6]
\draw (1,0) node {} -- (1.5,0) ;\draw (0.5,0.86602540378444) node{} -- (0.75,1.2990381056767);\draw (2,0.86602540378444) node{} -- (1.75,1.2990381056767);\fill  (1,0) circle[radius=1.3pt]  (1.5,0)  circle[radius=1.3pt]  (0.5,0.86602540378444) circle[radius=1.3pt]  (0.75,1.2990381056767) circle[radius=1.3pt]   (2,0.86602540378444) circle[radius=1.3pt]  (1.75,1.2990381056767) circle[radius=1.3pt]     ;\draw  (0,0) node (11) {} -- (1,0) node (12) {} -- (0.5,0.86602540378444) node (13) {} node[pos=0.5] (edge12-13) {} --cycle;\fill  (0,0) circle[radius=1.3pt]  (1,0)  circle[radius=1.3pt]  (0.5,0.86602540378444) circle[radius=1.3pt] ;\draw  (1.5,0) node[coordinate] (21) {} -- (2.5,0) node (22) {} -- (2,0.86602540378444) node (23) {}-- (21) node[pos=0.5] (edge22-21) {}  ;\fill  (1.5,0) circle[radius=1.3pt]  (2.5,0)  circle[radius=1.3pt]  (2,0.86602540378444) circle[radius=1.3pt] ;\draw  (0.75,1.2990381056767) node (31) {} -- (1.75,1.2990381056767) node[pos=0.5] (edge31-32) {} node (32) {} -- (1.25,2.1650635094611) node (33) {}--cycle;\fill  (0.75,1.2990381056767) circle[radius=1.3pt]  (1.75,1.2990381056767)  circle[radius=1.3pt]  (1.25,2.1650635094611) circle[radius=1.3pt] ;
\draw (0,0) to  [out=180,in=120]  (-.433,-.25)  to [out=300, in=240] (0,0);
\draw (2.5,0) to [out=0,in=60]   (2.933,-.25)  to [out=240, in=300] (2.5,0);
\draw (1.25,2.16506) to  [out=60,in=0]  (1.25,2.66506)  to  [out=180, in=120] (1.25,2.16506);
\path (edge12-13) node[right]{$y'$}  (edge22-21) node[left]{$y'$}  (edge31-32) node[below]{$y'$}  ;
\end{tikzpicture}
\ \ \ \ 
\begin{tikzpicture}[scale=1.6]
\draw (1,0) node{} -- (1.5,0);\draw (0.5,0.86602540378444) node{} -- (0.75,1.2990381056767);\draw (2,0.86602540378444) node{} -- (1.75,1.2990381056767);\fill  (1,0) circle[radius=1.3pt]  (1.5,0)  circle[radius=1.3pt]  (0.5,0.86602540378444) circle[radius=1.3pt]  (0.75,1.2990381056767) circle[radius=1.3pt]   (2,0.86602540378444) circle[radius=1.3pt]  (1.75,1.2990381056767) circle[radius=1.3pt]     ;\draw (0.4,0) node{} -- (0.6,0);\draw (0.2,0.34641016151378) node{} -- (0.3,0.51961524227066);\draw (0.8,0.34641016151378) node{} -- (0.7,0.51961524227066);\fill  (0.4,0) circle[radius=1.3pt]  (0.6,0)  circle[radius=1.3pt]  (0.2,0.34641016151378) circle[radius=1.3pt]  (0.3,0.51961524227066) circle[radius=1.3pt]   (0.8,0.34641016151378) circle[radius=1.3pt]  (0.7,0.51961524227066) circle[radius=1.3pt]     ;\draw  (0,0) node (111) {} -- (0.4,0)  node (112) {} -- (0.2,0.34641016151378) node[pos=0.5] (edge112113) {} node (113) {}--cycle;\fill  (0,0) circle[radius=1.3pt]  (0.4,0)  circle[radius=1.3pt]  (0.2,0.34641016151378) circle[radius=1.3pt] ;\draw  (0.6,0) node[coordinate] (121) {} -- (1,0) node (122) {} -- (0.8,0.34641016151378) node (123) {}-- (121) node[pos=0.5] (edge123121) {};\fill  (0.6,0) circle[radius=1.3pt]  (1,0)  circle[radius=1.3pt]  (0.8,0.34641016151378) circle[radius=1.3pt] ;\draw  (0.3,0.51961524227066) node (131) {} -- (0.7,0.51961524227066) node[pos=0.5] (edge131132) {} node (132) {} -- (0.5,0.86602540378444) node (133) {}--cycle;\fill  (0.3,0.51961524227066) circle[radius=1.3pt]  (0.7,0.51961524227066)  circle[radius=1.3pt]  (0.5,0.86602540378444) circle[radius=1.3pt] ;\draw (1.9,0) node{} -- (2.1,0);\draw (1.7,0.34641016151378) node{} -- (1.8,0.51961524227066);\draw (2.3,0.34641016151378) node{} -- (2.2,0.51961524227066);\fill  (1.9,0) circle[radius=1.3pt]  (2.1,0)  circle[radius=1.3pt]  (1.7,0.34641016151378) circle[radius=1.3pt]  (1.8,0.51961524227066) circle[radius=1.3pt]   (2.3,0.34641016151378) circle[radius=1.3pt]  (2.2,0.51961524227066) circle[radius=1.3pt]     ;\draw  (1.5,0) node (211) {} -- (1.9,0) node (212) {} -- (1.7,0.34641016151378)  node[pos=0.5] (edge212213) {} node (213) {}--cycle;\fill  (1.5,0) circle[radius=1.3pt]  (1.9,0)  circle[radius=1.3pt]  (1.7,0.34641016151378) circle[radius=1.3pt] ;\draw  (2.1,0) node (221) {} -- (2.5,0) node (222) {} -- (2.3,0.34641016151378) node (223) {}--cycle node[pos=0.5] (edge223221) {} ;\fill  (2.1,0) circle[radius=1.3pt]  (2.5,0)  circle[radius=1.3pt]  (2.3,0.34641016151378) circle[radius=1.3pt] ;\draw  (1.8,0.51961524227066) node (231) {} -- (2.2,0.51961524227066) node[pos=0.5] (edge231232) {} node (232) {} -- (2,0.86602540378444) node (233) {}--cycle;\fill  (1.8,0.51961524227066) circle[radius=1.3pt]  (2.2,0.51961524227066)  circle[radius=1.3pt]  (2,0.86602540378444) circle[radius=1.3pt] ;\draw (1.15,1.2990381056767) node{} -- (1.35,1.2990381056767);\draw (0.95,1.6454482671904) node{} -- (1.05,1.8186533479473);\draw (1.55,1.6454482671904) node{} -- (1.45,1.8186533479473);\fill  (1.15,1.2990381056767) circle[radius=1.3pt]  (1.35,1.2990381056767)  circle[radius=1.3pt]  (0.95,1.6454482671904) circle[radius=1.3pt]  (1.05,1.8186533479473) circle[radius=1.3pt]   (1.55,1.6454482671904) circle[radius=1.3pt]  (1.45,1.8186533479473) circle[radius=1.3pt]     ;\draw  (0.75,1.2990381056767) node (311) {} -- (1.15,1.2990381056767) node (312) {} -- (0.95,1.6454482671904) node[pos=0.5] (edge312313) {}  node (313) {}--cycle;\fill  (0.75,1.2990381056767) circle[radius=1.3pt]  (1.15,1.2990381056767)  circle[radius=1.3pt]  (0.95,1.6454482671904) circle[radius=1.3pt] ;\draw  (1.35,1.2990381056767) node (321) {} -- (1.75,1.2990381056767) node (322) {} -- (1.55,1.6454482671904) node (323) {}--cycle node[pos=0.5] (edge323321) {};\fill  (1.35,1.2990381056767) circle[radius=1.3pt]  (1.75,1.2990381056767)  circle[radius=1.3pt]  (1.55,1.6454482671904) circle[radius=1.3pt] ;\draw  (1.05,1.8186533479473) node (331) {} -- (1.45,1.8186533479473)node[pos=0.5] (edge331332) {}  node (332) {} -- (1.25,2.1650635094611) node (333) {}--cycle;\fill  (1.05,1.8186533479473) circle[radius=1.3pt]  (1.45,1.8186533479473)  circle[radius=1.3pt]  (1.25,2.1650635094611) circle[radius=1.3pt] ;
\draw (0,0) to  [out=180,in=120]  (-.433,-.25)  to [out=300, in=240] (0,0);
\draw (2.5,0) to [out=0,in=60]   (2.933,-.25)  to [out=240, in=300] (2.5,0);
\draw (1.25,2.16506) to  [out=60,in=0]  (1.25,2.66506)  to  [out=180, in=120] (1.25,2.16506);
\path (edge112113) node[right]  {$\scriptstyle y$}  (edge123121) node[left]  {$\scriptstyle y$}  (edge131132) node[below] {$\scriptstyle y$}   (edge212213) node[right]  {$\scriptstyle y$}  (edge223221) node[left]  {$\scriptstyle y$}  (edge231232) node[below] {$\scriptstyle y$} (edge312313) node[right]  {$\scriptstyle y$}  (edge323321) node[left]  {$\scriptstyle y$}  (edge331332) node[below] {$\scriptstyle y$}  ;    
\end{tikzpicture}
\caption{Weighted graphs on $H_2$ and $H_3$; unmarked edges have weight $1$.}\label{fig:GSmethod}
\end{figure}

Unfortunately, direct computation immediately shows that this Schur complement does not provide a spectral similarity between  $\Delta^H_{n+1}$ and $\Delta^H_n$; indeed, the Schur complement corresponds to the Laplacian on $H_n$ for a weighted graph in which the edges between pairs of vertices both of which are in $\tilde V^H_n$ are different than the weight on all other edges. This motivates the introduction of such a weighted Laplacian, which we illustrate in Figure~\ref{fig:GSmethod}.  Grigorchuk and 
\v{S}uni\'{c} then show this weighted Laplacian on $H_{n+1}$ is spectrally similar to the weighted Laplacian on $H_n$, only now the decimation function depends on both the spectral value $x$ and the weight $y$. The remarkable thing is that they are then able to find an explicit semiconjugacy of this two dimensional dynamics for the weighted Laplacian to the one-dimensional dynamics we found in Theorem~\ref{thm:Hndecimation}.

This brings us to the difference between our methods, which is that our vertex set $V^H_n$ does not lie in $V^H_{n+1}$, and we therefore are not  taking the Schur complement with respect to the subspace of functions supported on $\tilde V^H_{n+1}\setminus \tilde V^H_{n}$. Indeed, the subspace with respect to which we are taking the Schur complement does not correspond to any subset of the vertices of the graph $H_{n+1}$.  As far as the authors are aware, this is the first example using such a spectral decimation on self-similar graphs, though the use of spectral decimation with respect to subspaces coming from vertex subsets has been in common use for many years.  Given the increased flexibility obtained, it seems this variation of the spectral decimation method warrants further investigation.

\section{Acknowledgements}
The authors thank Alexander Teplyaev for suggesting this enjoyable problem, the University of Connecticut for hosting the authors, and grants NSF DMS-1659643 and DMS-1613025 for financial support that made this work possible.


\begin{thebibliography}{10}

\bibitem{BajorinChenJPhs}
N.~Bajorin, T.~Chen, A.~Dagan, C.~Emmons, M.~Hussein, M.~Khalil, P.~Mody,
  B.~Steinhurst, and A.~Teplyaev.
\newblock Vibration modes of {$3n$}-gaskets and other fractals.
\newblock {\em J. Phys. A}, 41(1):015101, 21, 2008.

\bibitem{BG}
L.~Bartholdi and R.~I. Grigorchuk.
\newblock On the spectrum of {H}ecke type operators related to some fractal
  groups.
\newblock {\em Tr. Mat. Inst. Steklova}, 231(Din. Sist., Avtom. i Beskon.
  Gruppy):5--45, 2000.

\bibitem{DSV}
Kyallee Dalrymple, Robert~S. Strichartz, and Jade~P. Vinson.
\newblock Fractal differential equations on the {S}ierpinski gasket.
\newblock {\em J. Fourier Anal. Appl.}, 5(2-3):203--284, 1999.

\bibitem{FukushimaShima}
M.~Fukushima and T.~Shima.
\newblock On a spectral analysis for the {S}ierpi\'nski gasket.
\newblock {\em Potential Anal.}, 1(1):1--35, 1992.

\bibitem{GrigorchukSunic}
Rostislav Grigorchuk and Zoran \v{S}uni\'{c}.
\newblock Schreier spectrum of the {H}anoi {T}owers group on three pegs.
\newblock In {\em Analysis on graphs and its applications}, volume~77 of {\em
  Proc. Sympos. Pure Math.}, pages 183--198. Amer. Math. Soc., Providence, RI,
  2008.

\bibitem{KSW}
Daniel~J. Kelleher, Benjamin~A. Steinhurst, and Chuen-Ming~M. Wong.
\newblock From self-similar structures to self-similar groups.
\newblock {\em Internat. J. Algebra Comput.}, 22(7):1250056, 16, 2012.

\bibitem{Kigbook}
Jun Kigami.
\newblock {\em Analysis on fractals}, volume 143 of {\em Cambridge Tracts in
  Mathematics}.
\newblock Cambridge University Press, Cambridge, 2001.

\bibitem{MT}
Leonid Malozemov and Alexander Teplyaev.
\newblock Self-similarity, operators and dynamics.
\newblock {\em Math. Phys. Anal. Geom.}, 6(3):201--218, 2003.

\bibitem{Nekbook}
Volodymyr Nekrashevych.
\newblock {\em Self-similar groups}, volume 117 of {\em Mathematical Surveys
  and Monographs}.
\newblock American Mathematical Society, Providence, RI, 2005.

\bibitem{NekTep}
Volodymyr Nekrashevych and Alexander Teplyaev.
\newblock Groups and analysis on fractals.
\newblock In {\em Analysis on graphs and its applications}, volume~77 of {\em
  Proc. Sympos. Pure Math.}, pages 143--180. Amer. Math. Soc., Providence, RI,
  2008.

\bibitem{RammalToulouse}
R.~Rammal and G.~Toulouse.
\newblock Random walks on fractal structures and percolation clusters.
\newblock {\em J. Phys. Lett.}, 44:L13–L22, 1983.

\bibitem{StrichartzAverages}
Robert~S. Strichartz.
\newblock The {L}aplacian on the {S}ierpinski gasket via the method of
  averages.
\newblock {\em Pacific J. Math.}, 201(1):241--256, 2001.

\bibitem{StrichartzFractafolds}
Robert~S. Strichartz.
\newblock Fractafolds based on the {S}ierpi\'{n}ski gasket and their spectra.
\newblock {\em Trans. Amer. Math. Soc.}, 355(10):4019--4043, 2003.

\bibitem{Strichartzbook}
Robert~S Strichartz.
\newblock {\em Differential equations on fractals: a tutorial}.
\newblock Princeton University Press, 2006.

\bibitem{StrichartzTransformations}
Robert~S. Strichartz.
\newblock Transformation of spectra of graph {L}aplacians.
\newblock {\em Rocky Mountain J. Math.}, 40(6):2037--2062, 2010.

\bibitem{StrichartzTeplyaev}
Robert~S. Strichartz and Alexander Teplyaev.
\newblock Spectral analysis on infinite {S}ierpi\'{n}ski fractafolds.
\newblock {\em J. Anal. Math.}, 116:255--297, 2012.

\end{thebibliography}
\end{document}